# EDGE EXCHANGEABLE MODELS FOR NETWORK DATA

HARRY CRANE AND WALTER DEMPSEY

ABSTRACT. Exchangeable models for countable vertex-labeled graphs cannot replicate the large sample behaviors of sparsity and power law degree distribution observed in many network datasets. Out of this mathematical impossibility emerges the question of how network data can be modeled in a way that reflects known empirical behaviors and respects basic statistical principles. We address this question by observing that edges, not vertices, act as the statistical units in networks constructed from interaction data, making a theory of edge-labeled networks more natural for many applications. In this context we introduce the concept of *edge exchangeability*, which unlike its vertex exchangeable counterpart admits models for networks with sparse and/or power law structure. Our characterization of edge exchangeable networks gives rise to a class of nonparametric models, akin to graphon models in the vertex exchangeable setting. Within this class, we identify a tractable family of distributions with a clear interpretation and suitable theoretical properties, whose significance in estimation, prediction, and testing we demonstrate.

## 1. INTRODUCTION

Statistical network analysis is hamstrung by the lack of an inferential framework that both facilitates sound conclusions from and replicates empirical properties of network data. Of the basic challenges facing network analysis, Robins and Morris (2007) wrote,

> *A good [network] model needs to be both estimable from data and a reasonable representation of that data[...]. Models that cannot be estimated from the data, that cannot reproduce the data to some adequate extent, models with implausible or empirically unsound assumptions, [...] cannot be considered complete.* (Robins and Morris, 2007, p. 169)

Above all, this quote emphasizes that sound modeling precedes, and is a crucial part of, methodological progress in the study of networks. The importance of modeling is further underscored by the fact that much of the current body of statistical theory and methods for networks has been developed for a select class of tractable models, such as graphon models

*Date*: October 21, 2016.

H. Crane is partially supported by NSF grants CNS-1523785 and CAREER DMS-1554092.

R code and links to datasets publicly available at https://github.com/wdempsey/e2.

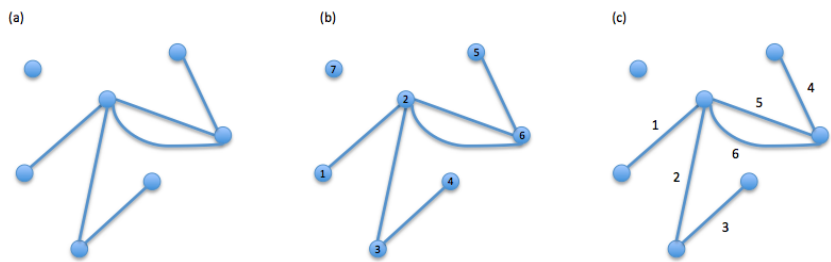

Figure 1. (a) A network structure derived from some process of interactions. Neither vertices nor edges come equipped with labels. (b) Network obtained by labeling vertices of network in Panel (a). (c) Network obtained by labeling edges of network in Panel (a).

(Bickel and Chen, 2009; Gao et al., 2015; Klopp et al., 2016; Wolfe and Olhede, 2014) and the stochastic blockmodel (Gao et al., 2016; Zhao et al., 2011), which were not originally intended to handle complex networks and tend to fit basic features of modern datasets rather poorly.

Here we initiate a framework for modeling and analyzing network data constructed from interaction processes, such as communications (Klimt and Yang, 2004; Opsahl and Panzarasa, 2009), collaborations (Barabási and Albert, 1999; Newman, 2001), and relationships (McAuley and Leskovec, 2012; Zachary, 1977). In this context, we assume the relevant information for the interactions is captured by an unlabeled network structure, as in Figure 1(a), with the vertices corresponding to elements of a population and edges representing interactions among adjacent vertices.

The convention in statistics, and network science more broadly, is to represent the structure in Figure 1(a) as a graph with labeled vertices, as in Figure 1(b). Out of this convention emerges the tendency—on grounds of logic as well as expediency—to assume an exchangeable model, which assigns equal probability to any two graphs that are equivalent up to relabeling vertices; see Figure 2.

On one hand, exchangeability seems to correct for any arbitrariness in labeling the vertices as well as narrow the analysis to a tractable class of models. On the other hand, exchangeable random graph models cannot replicate the widely observed empirical properties of sparsity and power law degree distribution observed in many modern network datasets (Abello et al., 1998; Barabási and Albert, 1999; Faloutsos et al., 1999; Jeong et al., 2001). Paraphrasing a well known outcome from probability and combinatorics (Diaconis and Janson, 2008; Lovász and Szegedy, 2006):

> *A realization from an exchangeable model for countable graphs is sparse if and only if it has no edges.*

Together, these observations prompt a foundational question of network modeling:

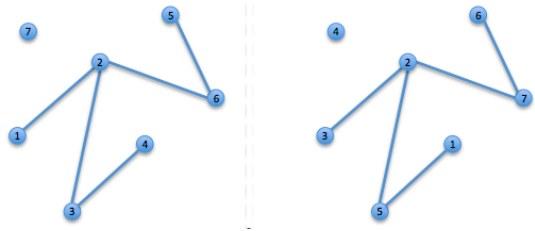

FIGURE 2. Two different ways to label vertices of the same unlabeled network. Vertex exchangeable models assign equal probability to both outcomes.

*Is there a notion of probabilistic symmetry whose ergodic measures [...] describe useful statistical models for sparse graphs with network properties?* (Orbanz and Roy, 2015, p. 459)

We address this question as part of our broader development of *edge exchangeable network models*, which are most appropriate when the edges are the statistical units, as they are for the interaction processes we study. We define this setup more precisely in Sections 2-3 and go on to establish basic properties of edge exchangeable models throughout Sections 4-8.

## 2. Interaction data

**Definition 2.1** (Interaction data). *For a set $\mathcal{P}$, we write fin($\mathcal{P}$) to denote the set of all finite (ordered) multisets of $\mathcal{P}$. An* interaction process *for a population $\mathcal{P}$ is a correspondence $\mathcal{I}: I \to$ fin($\mathcal{P}$) between a set $I$ indexing interactions and finite multisets of $\mathcal{P}$.*

The interaction processes in Definition 2.1 are the primary objects of our study, for which we assume no additional information, such as covariates. Incorporating covariates is notoriously difficult in network analysis, with only a limited number of successful attempts to date (Airoldi et al, 2011; Latouche et al, 2015; Mariadassou et al, 2010; Sweet, 2015; Tallberg, 2004; Zhang et al, 2015). Butts (2008) and Perry and Wolfe (2013) discuss elements of social theory in the context of interaction networks.

In Definition 2.1, each $\mathcal{I}(i)$ consists of the individuals involved in interaction $i \in I$. The ordering within $\mathcal{I}(i)$ captures directional interactions present in some applications, for example, the sending of emails from one account to a list of recipients. This setup contains the case in which interactions are undirected, which we obtain by specifying a canonical ordering of the elements in $\mathcal{P}$ and listing interactions according to that ordering.

Definition 2.1 captures the structure of many interaction datasets studied in the networks literature. In the actors collaboration network (Barabási and Albert, 1999; Rossi and Ahmed, 2015), for example, $I$ indexes a set of movies, with each movie $i \in I$ corresponding to an interaction involving the set of actors $\mathcal{I}(i)$ in its cast. In the Enron network (Klimt and Yang, 2004), interactions correspond to emails exchanged among Enron employees. Crane and Dempsey (2016, Section B) list other common interaction datasets.

2.1. **Network representation of interaction data.** Though not inherently structured as networks, interaction datasets are often depicted as in Figure 1. An interaction process $\mathcal{I} : I \to \text{fin}(\mathcal{P})$ is commonly represented by a *graph* $G_{\mathcal{I}} = (\mathcal{P}, E)$ with vertex set $\mathcal{P}$ and edge set $E \subseteq \mathcal{P} \times \mathcal{P}$ defined by

(1) $\quad (v, v') \in E \quad \text{if and only if} \quad \{v, v'\} \subseteq \mathcal{I}(i) \text{ for some } i \in I;$

that is, there is an edge between $v$ and $v'$ in $G_{\mathcal{I}}$ if $v$ and $v'$ both participate in at least one common interaction in $\mathcal{I}$.

The representation in (1) may disregard important features of the data, such as interactions involving more than two individuals and repeated interactions among the same set of individuals. In the actors dataset, for example, there are often more than two actors in each movie, nothing precludes the same set of actors from being cast together in more than one movie, and actors sometimes play more than one role in the same film.

These features of $\mathcal{I} : I \to \text{fin}(\mathcal{P})$ could be more faithfully represented by a hypergraph $H_{\mathcal{I}} : \text{fin}(\mathcal{P}) \to \{0, 1, \ldots\}$, where

(2) $\quad H_{\mathcal{I}}(A) = \#\{i \in I : \mathcal{I}(i) = A\}, \quad A \in \text{fin}(\mathcal{P}),$

records the number of interactions corresponding to each multiset $A$. This representation, however, presents another practical problem, as most statistical models are tailored specifically to graphs, with little or no easy extension to handle hyperedges with multiplicity. See Kivelä et al. (2014) for some discussion of hypergraph models in the context of multilayer networks.

2.2. **Statistical units.** The convention of representing network data as a graph or hypergraph with labeled vertices implicitly, and in most cases unintentionally, treats the vertices as the statistical units for the intended application. To be clear, the mere assignment of labels to the observed vertices, as in Figure 1(b), presents no immediate issue. The trouble arises only when the data are assumed to be a partial observation of a larger population process, as is often the case.

In this situation, the units are inherently tied to sampling, in that they represent the smallest elements on which observations are taken. The act of labeling vertices and treating the observed data as a graph $G = (V, E)$, as in (1), is consistent with the assumption that $G$ is a partial observation of a population graph $G' = (V', E')$ with $V' \supseteq V$ and $E' \subseteq V' \times V'$, so that the data $G = G'|_V$ corresponds to the restriction of $G'$ to vertices $V$. In this way, the observed data $G$, even if not explicitly stated or intended, is assumed to contain information about all the edges between sampled vertices.

For interaction data, however, the interactions, and therefore edges, comprise the fundamental units of observation. For example, we may observe the actors dataset by sampling movies without replacement from the Internet Movie Database[1] (IMDB). Or, in the Enron dataset, we may observe only those emails exchanged during a specific period of time, which constitutes a sample from the collection of all emails exchanged in the company's history. In either case, the data result from a partial observation of the interactions.

For this purpose, we assume a population process $\mathcal{I} : \mathbb{N} \to \text{fin}(\mathcal{P})$ of a potentially infinite collection of interactions indexed by the positive integers $\mathbb{N} = \{1, 2, \ldots\}$. From $\mathcal{I}$, we observe the restriction $\mathcal{I}|_S : S \to \text{fin}(\mathcal{P})$, $s \mapsto \mathcal{I}(s)$, to some finite sample $S \subset \mathbb{N}$. The apparent role of interactions, and therefore edges, as the statistical units suggests the alternative representation as an edge-labeled network, as in Figure 1(c).

2.3. **Edge-labeled networks.** The foregoing discussion stresses the difference between the *population process*, as a collection of interactions $\mathcal{I} : \mathbb{N} \to \text{fin}(\mathcal{P})$, and its representation as so-called *network data* in Figure 1.

Let $\mathcal{I} : \mathbb{N} \to \text{fin}(\mathcal{P})$ be an interaction process. Any bijection $\rho : \mathcal{P} \to \mathcal{P}'$ induces an action $\text{fin}(\mathcal{P}) \to \text{fin}(\mathcal{P}')$ by the map

$$(3) \qquad s = (s_1, \ldots, s_r) \mapsto \rho s = (\rho(s_1), \ldots, \rho(s_r)) \in \text{fin}(\mathcal{P}').$$

With this, any bijection $\rho : \mathcal{P} \to \mathcal{P}'$ acts on $\mathcal{I} : \mathbb{N} \to \text{fin}(\mathcal{P})$ by composition of maps, $\mathcal{I} \mapsto \rho\mathcal{I}$, with $\rho\mathcal{I} : \mathbb{N} \to \text{fin}(\mathcal{P}')$ given by

$$(4) \qquad (\rho\mathcal{I})(i) = \rho(\mathcal{I}(i)), \quad i \in \mathbb{N},$$

as defined in (3). From this, we define the *edge-labeled network induced by* $\mathcal{I} : S \to \text{fin}(\mathcal{P})$ as the equivalence class

$$(5) \qquad \mathcal{E}_{\mathcal{I}} = \bigcup_{\mathcal{P}' : \#\mathcal{P}' = \#\mathcal{P}} \{\mathcal{I}' : S \to \text{fin}(\mathcal{P}') : \rho\mathcal{I} = \mathcal{I}' \text{ for some bijection } \rho : \mathcal{P} \to \mathcal{P}'\}.$$

The equivalence class $\mathcal{E}_{\mathcal{I}}$ corresponds to an edge-labeled network structure, as in Figure 1(c). Note that $\mathcal{E}_{\mathcal{I}}$ does not depend on the specific set $\mathcal{P}$, and so we may disregard $\mathcal{P}$, or treat it implicitly as $\mathcal{P} = \mathbb{N}$, in our discussion. We write $\mathfrak{E}_S$ for the set of networks with edges labeled in $S \subseteq \mathbb{N}$.

For any $S' \subseteq S$, we define the *restriction* of $\mathcal{E} \in \mathfrak{E}_S$ to $\mathfrak{E}_{S'}$ by $\mathcal{E}|_{S'}$, the edge-labeled network obtained by removing any edges labeled in $S \setminus S'$. If $\mathcal{E} = \mathcal{E}_{\mathcal{I}}$ for some interaction process $\mathcal{I} : S \to \text{fin}(\mathcal{P})$, then $\mathcal{E}|_{S'}$ is the edge-labeled network induced by the restricted process $\mathcal{I}|_{S'} : S' \to \text{fin}(\mathcal{P})$, $s \mapsto \mathcal{I}(s)$.

**Remark 2.2.** *For clarity we reserve the term* graph *to specifically refer to a vertex-labeled structure, such as the objects given in* (1), (2), *and Figure 1(b). We use the term* network *for the generic unlabeled structure in Figure 1(a) and* edge-labeled network *for the object defined in* (5) *and shown in Figure 1(c).*

---

[1] http://www.imdb.com/

3. Network properties

For any edge-labeled network $\mathcal{E}$, we define $v(\mathcal{E})$ as the number of non-isolated vertices in $\mathcal{E}$, that is, the number of vertices involved in at least one of the edges of $\mathcal{E}$. We also define $M_k(\mathcal{E})$ as the number of $k$-ary edges in $\mathcal{E}$ for each $k \geq 1$, $N_k(\mathcal{E})$ as the number of vertices that appear exactly $k$ times in $\mathcal{E}$, and $d(\mathcal{E}) = (d_k(\mathcal{E}))_{k \geq 1}$ as the degree distribution of $\mathcal{E}$, where $d_k(\mathcal{E}) = N_k(\mathcal{E})/v(\mathcal{E})$ is the proportion of vertices with degree $k$ in $\mathcal{E}$. For example, the edge-labeled network $\mathcal{E}$ in Figure 1(c) has $v(\mathcal{E}) = 6$, $e(\mathcal{E}) = 6$, $M_k(\mathcal{E}) = 6$ for $k = 2$, $M_k(\mathcal{E}) = 0$ for $k \neq 2$, and $d(\mathcal{E}) = (3/6, 1/6, 1/6, 1/6, 0, \ldots)$. Notice that these statistics do not depend on the edge labels assigned to $\mathcal{E}$.

**Definition 3.1** (Sparsity and power law degree distribution). *Let $(\mathcal{E}_n)_{n \geq 1}$ be a sequence of edge-labeled networks for which $e(\mathcal{E}_n) \to \infty$ as $n \to \infty$. The sequence $(\mathcal{E}_n)_{n \geq 1}$ is* sparse *if*

$$\limsup_{n \to \infty} \frac{e(\mathcal{E}_n)}{v(\mathcal{E}_n)^{m_\bullet(\mathcal{E}_n)}} = 0, \tag{6}$$

*where $m_\bullet(\mathcal{E}_n) = e(\mathcal{E}_n)^{-1} \sum_{k \geq 1} k M_k(\mathcal{E}_n)$ is the average arity of the edges in $\mathcal{E}_n$.*

*The sequence $(\mathcal{E}_n)_{n \geq 1}$ exhibits a* power law degree distribution *if for some $\gamma > 1$ the degree distributions $(d(\mathcal{E}_n))_{n \geq 1}$ satisfy $d_k(\mathcal{E}_n) \sim \ell(k) k^{-\gamma}$ as $n \to \infty$ for all large $k$ for some slowly varying function $\ell(x)$, that is, $\lim_{x \to \infty} \ell(tx)/\ell(x) = 1$ for all $t > 0$, where $a_n \sim b_n$ indicates that $a_n/b_n \to 1$ as $n \to \infty$.*

**Remark 3.2.** *The slowly varying function $\ell(k)$ makes the definition of power law distribution robust to finite sample behavior and only affects the shape of the distribution. It does not alter the tail behavior.*

Our definition of sparsity in (6) refines the usual notion of sparsity for vertex-labeled graphs, as defined, for example, by Nesetril and Ossona de Mendez (2012). In the more familiar case when every interaction involves exactly 2 vertices, $m_\bullet(\mathcal{E}_n) = 2$ for all $n \geq 1$ and $(\mathcal{E}_n)_{n \geq 1}$ is sparse if

$$\limsup_{n \to \infty} \frac{e(\mathcal{E}_n)}{v(\mathcal{E}_n)^2} = 0.$$

4. Edge exchangeable network models

Given $\mathcal{E} \in \mathfrak{E}_S$, for $S \subseteq \mathbb{N}$, and a permutation $\sigma : S \to S$, we write $\mathcal{E}^\sigma$ to denote the edge-labeled network obtained by relabeling the edges of $\mathcal{E}$ according to $\sigma$. Figure 3 demonstrates this operation visually.

More precisely, given an edge-labeled network $\mathcal{E}$ based on interaction data $\mathcal{I} : S \to \text{fin}(\mathcal{P})$ and a permutation $\sigma : S \to S$, $\mathcal{E}^\sigma$ is the edge-labeled network induced by $\mathcal{I}^\sigma : S \to \text{fin}(\mathcal{P})$ with $\mathcal{I}^\sigma(i) = \mathcal{I}(\sigma^{-1}(i))$ for each $i \in S$. (It is clear from (5) that the definition of $\mathcal{E}^\sigma$ does not depend on which representative of the equivalence class $\mathcal{E}_\mathcal{I}$ is chosen when defining $\mathcal{E}^\sigma$.) Note the distinction between the action $\rho \mathcal{I}$ defined in (4) and $\mathcal{I}^\sigma$ given here. Whereas the action $\mathcal{I} \mapsto \rho \mathcal{I}$ in (4) reassigns labels to the population through a bijection $\rho : \mathcal{P} \to \mathcal{P}'$, the operation $\mathcal{I} \mapsto \mathcal{I}^\sigma$ corresponds to relabeling interactions and $\mathcal{E} \mapsto \mathcal{E}^\sigma$ to relabeling edges.

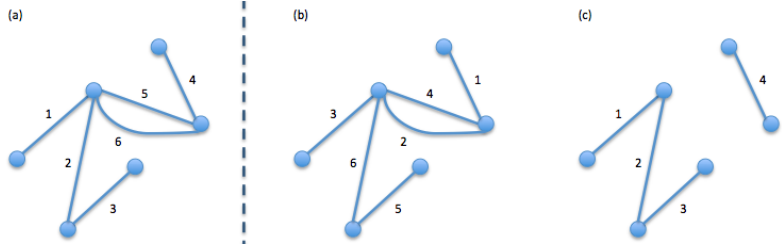

FIGURE 3. Relabeling and restriction operations for edge-labeled networks. Panel (a) shows a network with edges labeled $1, \ldots, 6$. Panel (b) shows the network from (a) with edges relabeled according to permutation $(1354)(26)$. Panel (c) shows the restriction of the network in (a) to the edges labeled in $1, \ldots, 4$.

We often write $\mathcal{Y}$ to denote a random edge-labeled network.

**Definition 4.1** (Edge exchangeable network). *A random edge-labeled network $\mathcal{Y} \in \mathfrak{E}_S$ is edge exchangeable if $\mathcal{Y}^\sigma =_\mathcal{D} \mathcal{Y}$ for all permutations $\sigma : S \to S$, where $=_\mathcal{D}$ denotes equality in distribution.*

Edge exchangeable models assign the same probability to all edge-labeled networks that are isomorphic up to relabeling, as shown in Figure 4.

The assumption of edge exchangeability has a natural interpretation in terms of sampling in that the interactions that determine an observed edge-labeled network $\mathcal{E} \in \mathfrak{E}_S$ are assumed to be representative of a larger population of interactions $\mathcal{I} : \mathbb{N} \to \text{fin}(\mathcal{P})$. Contrast this with the assumption of exchangeability for vertex-labeled graphs, in which case sampled vertices are assumed to be a representative sample of all vertices.

4.1. **Characterization of edge exchangeable networks.** We describe here a special case a more general de Finetti-type representation theorem of edge exchangeable networks, which we defer to Crane and Dempsey (2016, Theorem A.2). The case below is most relevant to statistical applications.

To construct an edge exchangeable random network $\mathcal{Y} \in \mathfrak{E}_\mathbb{N}$, we first choose any $f = (f_s)_{s \in \text{fin}(\mathbb{N})}$ in the *fin($\mathbb{N}$)-simplex*,

$$\mathcal{F}_1 = \left\{ (f_s)_{s \in \text{fin}(\mathbb{N})} : f_s \geq 0 \quad \text{and} \quad \sum_{s \in \text{fin}(\mathbb{N})} f_s = 1 \right\}.$$

Each $f \in \mathcal{F}_1$ determines a distribution on finite multisets of $\mathbb{N}$ by

(7) $$\text{pr}(X = s \mid f) = f_s, \quad s \in \text{fin}(\mathbb{N}),$$

which in turn determines an edge exchangeable network in $\mathfrak{E}_\mathbb{N}$ as follows.

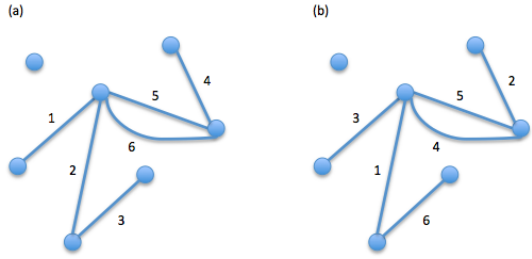

FIGURE 4. Two edge labelings of the network data from Figure 1. An edge exchangeable model assigns equal probability to both networks. Vertices are not labeled and, therefore, vertex labels play no role in the exchangeability condition.

Let $X_1, X_2, \ldots$ be independent, identically distributed (i.i.d.) random multisets drawn from (7). Given $X_1, X_2, \ldots$, define the random interaction process $\mathcal{X} : \mathbb{N} \to \text{fin}(\mathbb{N})$ by $\mathcal{X}(i) = X_i$ and write $\mathcal{E}_\mathcal{X}$ to denote the corresponding edge-labeled network obtained through (5), whose distribution we denote by $\epsilon_f$. See Figure 5(a) for an illustration.

With probability 1, each $s \in \text{fin}(\mathbb{N})$ occurs either zero or infinitely many times in the interaction process $\mathcal{X} : \mathbb{N} \to \text{fin}(\mathbb{N})$ corresponding to $X_1, X_2, \ldots$ chosen i.i.d. from (7), and the resulting edge-labeled network $\mathcal{E}_\mathcal{X}$ is edge exchangeable. This is a special case of the more generic edge exchangeable processes, for which interactions can occur zero, one, or infinitely many times. We call an interaction that occurs only once a *blip*, and we call an interaction process and its induced edge-labeled network *blip-free* if it contains no blips. We state here the blip-free version of our general representation theorem (Crane and Dempsey, 2016, Theorem A.2).

**Theorem 4.2.** *Let $\mathcal{Y} \in \mathfrak{E}_\mathbb{N}$ be an edge exchangeable network that is blip-free with probability 1. Then there exists a probability measure $\phi$ on $\mathcal{F}_1$ such that $\mathcal{Y} \sim \epsilon_\phi$, where*

$$\epsilon_\phi(\cdot) = \int_{\mathcal{F}_1} \epsilon_f(\cdot) \phi(df). \tag{8}$$

*That is, every blip-free edge exchangeable network $\mathcal{Y} \in \mathfrak{E}_\mathbb{N}$ can be generated by first sampling $f \sim \phi$ and, given $f$, putting $\mathcal{Y} = \mathcal{E}_\mathcal{X}$ with $\mathcal{X} : \mathbb{N} \to \text{fin}(\mathbb{N})$ defined by $\mathcal{X}(i) = X_i$ for $X_1, X_2, \ldots$ i.i.d. from (7).*

**Remark 4.3.** *The measure $\phi$ in Theorem 4.2 is not unique. Uniqueness results from a more technical treatment; see Crane and Dempsey (2016, Theorem A.2).*

**Corollary 4.4.** *The vertices in any blip-free edge exchangeable network $\mathcal{Y} \in \mathfrak{E}_\mathbb{N}$ arrive in size-biased order weighted by the relative frequency of their occurrence.*

According to Corollary 4.4, the assumption of edge exchangeability is incompatible with vertex exchangeability, under which the vertices may be interpreted as arriving in exchangeable random order. Edge exchangeable networks, therefore, do not admit a graphon representation in the sense of Lovász and Szegedy (2006). But the generic construction of edge exchangeable networks in (7) and Theorem 4.2 suggests the following analog to the graphon theory in terms of elements $f \in \mathcal{F}_1$.

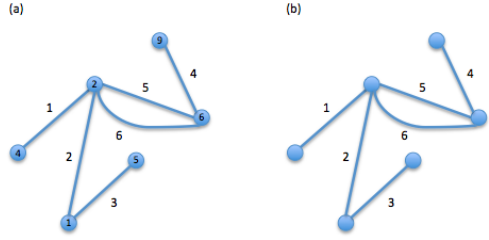

FIGURE 5. Illustration of the generic generating process for edge exchangeable networks described in (7) and Theorem 4.2. (a) Network representation of the interaction process $X_1 = \{2, 4\}$, $X_2 = \{1, 2\}$, $X_3 = \{1, 5\}$, $X_4 = \{6, 9\}$, $X_5 = \{2, 6\}$, $X_6 = \{2, 6\}$. (b) Edge-labeled network obtained by removing vertex labels from the network in Panel (a).

4.2. **Nonparametric edge exchangeable models.** Here we highlight the tractable nonparametric subclass of *vertex components models*, for which we construct $f = (f_s)_{s \in \text{fin}(\mathbb{N})}$ by specifying a probability distribution $\nu = \{\nu_k\}_{k \geq 1}$ on the positive integers, choosing a sequence $W = (W_i)_{i \geq 1}$ at random from the infinite simplex

$$\Delta_1 = \left\{(s_1, s_2, \ldots) : \sum_{i \geq 1} s_i = 1\right\},$$

and then defining

(9) $$f_{(s_1,\ldots,s_k)} = \nu_k \prod_{j=1}^{k} W_{s_j}, \quad (s_1, \ldots, s_k) \in \text{fin}(\mathbb{N}).$$

Under (9), we interpret the distribution in (7) as first choosing the size of a random multiset $s$ according to $\nu$ and, given $\#s = k$, choosing the $k$ elements i.i.d. from $(W_i)_{i \geq 1}$. (Recall from Definition 2.1 that fin($\mathbb{N}$) consists of ordered multisets. We recover the unordered case by writing each $\{s_1, \ldots, s_k\}$ so that $s_1 \leq s_2 \leq \cdots \leq s_k$ and multiplying (9) by the number of distinct ways to order the elements of $s$, which can be computed explicitly.)

Stick-breaking representations, see, for example, Sethuraman (1994), offer a computationally tractable way to fit the vertex components model to data. In this construction, we can generate the sequence $X_1, X_2, \ldots$ of interactions simultaneously with the vertex components $W = (W_i)_{i \geq 1}$.

Let $\{\varphi_i\}_{i \geq 1}$ be a collection of probability densities on $[0, 1]$. We first choose the size $K_n \sim \nu$ of the $n$th directed interaction $X_n = (X_{n,1}, \ldots, X_{n,K_n})$ independently for each $n = 1, 2, \ldots$. We begin by choosing $K_1 \sim \nu$, putting $X_{1,1} = 1$, and sampling $W_1 \sim \varphi_1$. We continue inductively for each $n = 1, 2, \ldots$ as follows. Given $K_n$ and $k < K_n$, we define $X_{n+1,\leq k}$ as the set of all $X_1, \ldots, X_n$ and $X_{n+1,1}, \ldots, X_{n+1,k}$ up to the $k$th chosen element of the $(n+1)$st interaction. Given $X_{n+1,\leq k}$ and $W_1, \ldots, W_{V_n}$, where $V_n = \max(X_{n+1,\leq k})$ is the largest vertex label assigned among $X_{n+1,\leq k}$, we choose $X_{n+1,k+1}$, provided $K_{n+1} \geq k+1$, according to

$$\text{pr}(X_{n+1,k+1} = r \mid X_{n+1,\leq k}, W_1, \ldots, W_{V_n}) = \begin{cases} W_r, & r = 1, \ldots, V_n, \\ 1 - \sum_{j=1}^{V_n} W_j, & r = V_n + 1. \end{cases} \tag{10}$$

If $X_{n+1,k+1} = V_n + 1$, then we choose $W_{V_n+1} \sim \varphi_{n+1}(\cdot/(1 - \sum_{j=1}^{V_n} W_j))$. (The division by $1 - \sum_{j=1}^{V_n} W_j$ represents the normalization by the remaining length of the stick in the stick-breaking construction.) We continue to draw $X_{n+1,k+2}, \ldots, X_{n+1,K_{n+1}}$ as in (10).

From the sequence $X_1, X_2, \ldots$ constructed above, we define $\mathcal{X}_n : [n] \to \text{fin}(\mathbb{N})$ as interaction process constructed by putting $\mathcal{X}_n(i) = X_i$ and defining $\mathcal{Y}_n = \mathcal{E}_{\mathcal{X}_n}$ to be the edge-labeled network induced by $\mathcal{X}_n$ for each $n \in \mathbb{N}$. We compute the joint density of $\mathcal{Y}_n$ and $W$ by

$$\text{pr}(\mathcal{Y}_n = \mathcal{E}, (W_1, \ldots, W_{v(\mathcal{E})}) \in (dw_i)_{1 \leq i \leq v(\mathcal{E})}; \{\varphi_i\}_{i \geq 1}) = \tag{11}$$

$$= \prod_{j=1}^{v(\mathcal{E})} \left(1 - \sum_{i=1}^{j-1} w_i\right) \varphi_j\left(\frac{w_j}{1 - \sum_{i=1}^{j-1} w_i}\right) w_j^{D_n(j)-1} dw_1 \cdots dw_{v(\mathcal{E})},$$

where $D_n(j)$ is the number of times the vertex corresponding to weight $w_j$ appears in $\mathcal{E}$ and $\sum_{i=1}^{0} w_i = 0$.

Below we study a two parameter family of models corresponding to the vertex components model in (11) with $W$ generated from a Dirichlet process (Ferguson, 1973). We explain this connection further in Section 5.2.

## 5. The Hollywood model

Though not limited to any specific application, we phrase the following *Hollywood model* in the terms of the actors collaboration network to aid comprehension. In this description, vertices correspond to actors and each edge contains the set of actors involved in the corresponding movie.

Let $\nu = \{\nu_k\}_{k \geq 1}$ be a probability distribution on the positive integers and let $(\alpha, \theta)$ satisfy either

- $0 < \alpha < 1$ and $\theta > -\alpha$ (infinite population) or
- $\alpha < 0$ and $\theta = -k\alpha$ for some $k = 1, 2, \ldots$ (finite population).

We generate a sequence of interactions $X_1, X_2, \ldots$ as follows.

Given $X_1, \ldots, X_{n-1}$, we choose the number of roles in the next movie independently according to $K_n \sim \nu$. Given $K_n = k$, we choose the $k$ actors $X_{n,1}, \ldots, X_{n,k}$ in order of their prominence, first filling the lead role, then the second lead role, and so on until all $k$ roles are filled. Let $V_n(j)$ be the number of unique actors seen in all the movies through the $(j-1)$st actor cast in movie $n$. (Thus, $V_n(1)$ is the number of unique actors appearing in movies $1,\ldots,n-1$.) For each $i = 1,\ldots,V_n(j)$, we write $D_{n,j}(i)$ to denote the number of roles for which the actor labeled $i$ has been cast up to and including the $(j-1)$st role of movie $n$. (Note that an actor may play more than one role in a given movie.) The actor $X_{n,j}$ cast in the $j$th role of movie $n$ is chosen randomly among the actors labeled $1,\ldots,V_n(j)$ and a previously unseen actor, labeled $V_n(j)+1$, according to

$$\text{pr}(X_{n,j} = i \mid X_1,\ldots,X_{n-1},X_{n,1},\ldots,X_{n,j-1}) \propto \begin{cases} D_{n,j}(i) - \alpha, & i = 1,\ldots,V_n(j), \\ \theta + \alpha V_n(j), & i = V_n(j)+1. \end{cases} \quad (12)$$

We update according to (12) until all $k$ roles of movie $n$ have been filled.

Each interaction $X_1, X_2, \ldots$ records the order in which actors are selected in each movie. From $X_1, X_2, \ldots$, we define $\mathcal{X}_n : [n] \to \text{fin}(\mathbb{N})$ by $\mathcal{X}_n(i) = X_i$ for $i = 1, \ldots, n$, and we write $\mathcal{Y}_n = \mathcal{E}_{\mathcal{X}_n}$ for the edge-labeled network induced by $\mathcal{X}_n$ as defined in (5). We call the resulting sequence $(\mathcal{Y}_n)_{n\geq 1}$ the *Hollywood process* with parameter $(\alpha, \theta, \nu)$, which determines a family of distributions on $\mathfrak{E}_{\mathbb{N}}$, called the *Hollywood model*.

The Hollywood model has a closed form expression for random edge-labeled networks of each finite size $n \geq 1$ given by

$$\text{pr}(\mathcal{Y}_n = \mathcal{E}; \alpha, \theta, \nu) =$$
$$= \left[\prod_{k \geq 1} \nu_k^{M_k(\mathcal{E})}\right] \alpha^{v(\mathcal{E})} \frac{(\theta/\alpha)^{\uparrow v(\mathcal{E})}}{\theta^{\uparrow m(\mathcal{E})}} \prod_{k=2}^{\infty} \exp\{N_k(\mathcal{E}) \log((1-\alpha)^{\uparrow (k-1)})\}, \quad (13)$$

where $\mathcal{E} \in \mathfrak{E}_{[n]}$, $v(\mathcal{E})$ is the number of non-isolated vertices in $\mathcal{E}$, $(N_k(\mathcal{E}))_{k\geq 1}$ gives the number of vertices with each degree $k \geq 1$, $M_k(\mathcal{E})$ is the number of $k$-ary edges in $\mathcal{E}$, $m(\mathcal{E}) = \sum_{k\geq 1} k M_k(\mathcal{E})$ is the total degree of $\mathcal{E}$, and $x^{\uparrow j} = x(x+1)\cdots(x+j-1)$ is the ascending factorial function.

Though defined for edge-labeled networks with directed edges, the Hollywood model determines a distribution on undirected edge-labeled networks by ignoring the edge orientations. In this case, the probability of an undirected network $\mathcal{E}^*$ is given by the expression in (13) multiplied by a combinatorial factor $C(\mathcal{E}^*)$ that counts the number of directed networks corresponding to $\mathcal{E}^*$. While we have no closed form expression for $C(\mathcal{E}^*)$, the quantity plays no role in our inferences because it does not depend on the parameters $(\alpha, \theta, \nu)$, as we discuss further in Section 6.1.

5.1. **Interpretation of parameters and finite population model.** The split parameter space of the Hollywood model enables both bounded and unbounded population sizes. The region $0 < \alpha < 1$ and $\theta > -\alpha$ gives rise to a sequence $(\mathcal{Y}_n)_{n\geq 1}$ for which $v(\mathcal{Y}_n) \to \infty$ almost surely (a.s.) as $n \to \infty$, as is reasonable to assume for the actors, Enron, and Wikipedia networks. The Karate Club dataset, on the other hand, is known to have exactly thirty-four club members but no limit on the number of interactions between individuals. The range $\alpha < 0$ and $\theta = -k\alpha$ accommodates this case by describing an edge exchangeable sequence $(\mathcal{Y}_n)_{n\geq 1}$ for which $v(\mathcal{Y}_n) \to k$ a.s. as $n \to \infty$.

By (12), $\alpha > 0$ increases the probability of observing previously unseen vertices but decreases the probability of observing a vertex again after its initial occurrence. The range $\alpha < 0$ has the opposite effect. Thus, $\alpha$ values near 1 make it more likely that new edges involve previously unseen vertices, but less likely that previously seen vertices occur in future edges. On the other hand, $\alpha < 0$ corresponds to a finite population size, so that each newly observed vertex decreases the number of unseen vertices and increases the probability that future edges involve previously seen vertices. In the $0 < \alpha < 1$ regime, larger values of $\theta$ increase the probability of seeing previously unobserved vertices in new edges, but the effect of $\theta$ diminishes as $n \to \infty$. In Section 5.3, we see that $0 < \alpha < 1$ is directly related to the sparsity and power law behavior of the sequence $(\mathcal{Y}_n)_{n\geq 1}$.

5.2. **Connection to vertex components model.** Crane (2016b) previously noted a connection between the Hollywood model and the Ewens–Pitman two parameter family of distributions on set partitions (Ewens, 1972; Pitman, 2005). The two models coincide in the unary setting of the Hollywood model, that is, the $(\alpha, \theta, \nu)$ case when $\nu_1 = 1$. In this sense, the Hollywood model may be viewed as a natural refinement of the Ewens–Pitman distribution and Chinese restaurant process, which enjoys wide relevance throughout statistics, mathematics, and applied science (Crane, 2016a).

For $\alpha < 0$ and $\theta = -k\alpha$, the Hollywood model corresponds to the vertex components model with $W = (W_1, \ldots, W_k)$ chosen from the symmetric Dirichlet distribution with parameter $(\alpha, \ldots, \alpha)$ on the $(k-1)$-simplex.

For $0 < \alpha < 1$ and $\theta > -\alpha$, the Hollywood model is a special case of the vertex components model with $W = (W_i)_{i\geq 1}$ chosen from the Griffiths–Engen–McCloskey (GEM) distribution with parameter $(\alpha, \theta)$ on $\Delta_{\mathbf{1}}$; see Feng (2010), also Crane (2016a), for further details on the GEM distribution and its relationship to the Poisson–Dirichlet distribution. Alternatively, the Hollywood model with $0 < \alpha < 1$ and $\theta > -\alpha$ can be constructed by noting the stick-breaking construction of the GEM distribution and taking $\varphi_j$ to be the density of the Beta distribution with parameter $(1 - \alpha, \theta + j\alpha)$ for each $j \geq 1$ (Pitman, 2005). We recover (13) by marginalizing over $W$ in (11).

5.3. **Statistical properties of the Hollywood model.** The close connection between the Hollywood model and the Poisson–Dirichlet distribution through the vertex components interpretation of Section 5.2 addresses the main questions posed in Section 1. Theorems 5.2 and 5.3, in particular, follow from the corresponding power law behavior of the Ewens–Pitman process (Crane, 2016a; Pitman, 2005).

**Theorem 5.1.** *The Hollywood model with parameter $(\alpha, \theta, \nu)$ determines an edge exchangeable probability distribution on $\mathfrak{E}_\mathbb{N}$ for all $(\alpha, \theta, \nu)$ in the parameter space of the model.*

**Theorem 5.2.** *Let $(\mathcal{Y}_n)_{n\geq 1}$ obey the Hollywood process with parameter $(\alpha, \theta, \nu)$. For each $n \geq 1$, let $p_n(k) = N_k(\mathcal{Y}_n)/v(\mathcal{Y}_n)$, $k \geq 1$, be the empirical degree distribution of $\mathcal{Y}_n$, where $N_k(\mathcal{Y}_n)$ is the number of vertices with degree $k \geq 1$ and $v(\mathcal{Y}_n)$ is the number of non-isolated vertices in $\mathcal{Y}_n$. Then, for every $k \geq 1$,*

$$p_n(k) \sim \alpha k^{-(\alpha+1)}/\Gamma(1-\alpha) \quad a.s. \quad \text{as } n \to \infty,$$

*where $\Gamma(t) = \int_0^\infty x^{t-1} e^{-x} dx$ is the gamma function. That is, $(\mathcal{Y}_n)_{n\geq 1}$ has a power law degree distribution with exponent $\gamma = 1 + \alpha \in (1, 2)$.*

**Theorem 5.3.** *Let $(\mathcal{Y}_n)_{n\geq 1}$ obey the Hollywood process with parameter $(\alpha, \theta, \nu)$ for $0 < \alpha < 1$ and $\theta > -\alpha$. Then the expected number of vertices in $\mathcal{Y}_n$ satisfies*

$$(14) \qquad E(v(\mathcal{Y}_n)) \sim \frac{\Gamma(\theta+1)}{\alpha \Gamma(\theta+\alpha)} (\mu n)^\alpha \quad \text{as } n \to \infty,$$

*where $\mu = \sum_{k\geq 1} k \nu_k$ is the mean edge arity. Furthermore, if $1/\mu < \alpha < 1$, then $(\mathcal{Y}_n)_{n\geq 1}$ is almost surely sparse in the sense of* (6).

The power law behavior in the range $1 < \gamma < 2$ of the Hollywood model complements that of the preferential attachment model (Chung and Lu, 2006). Some authors (Barabási and Albert, 1999) suggest that $\gamma > 2$ is more prevalent in datasets that exhibit power law, but more recent work by Crane and Dempsey (2015) demonstrates empirically that the range $1 < \gamma < 2$ of the Hollywood model is common for interaction networks.

For a more general class of edge exchangeable models that admits sparsity and power law, we define the *$\alpha$-diversity*, $0 < \alpha < 1$, of an edge exchangeable process $(\mathcal{Y}_n)_{n\geq 1}$ by

$$(15) \qquad S_\alpha = \lim_{n\to\infty} v(\mathcal{Y}_n)/e(\mathcal{Y}_n)^\alpha,$$

provided the limit exists and is strictly positive and finite almost surely. Following Pitman (2005, Lemma 3.11), an edge exchangeable process $(\mathcal{Y}_n)_{n\geq 1}$ from the vertex components model constructed from $\nu = \{\nu_k\}_{k\geq 1}$ and $W = (W_i)_{i\geq 1}$ has $\alpha$-diversity if and only if the reordered sequence $W_{(1)} \geq W_{(2)} \geq \cdots \geq 0$ satisfies $W_{(i)} \sim Z i^{-1/\alpha}$ a.s. as $i \to \infty$ for a random variable $0 < Z < \infty$. In this case, $(\mathcal{Y}_n)_{n\geq 1}$ is sparse for $1/\mu < \alpha < 1$, with $\mu = \sum_{k\geq 1} k \nu_k$, and has power law degree distribution with $\gamma = \alpha + 1$.

5.4. **Projecting to a network without multiple edges.** Though multiple edges occur naturally in networks constructed from interaction processes, most network models are designed to handle only simple graphs, and many network datasets record only a single edge to indicate the occurrence of some positive number of interactions, as in (1). The fact that many network datasets are obtained by thresholding edge multiplicities is often ignored during data analysis; and the significance of this action on inference is underappreciated in the broader literature.

We define the $(t, c)$-*projection* of an interaction process $\mathcal{I} : \mathbb{N} \to \text{fin}(\mathcal{P})$ by $H_{\mathcal{I}}^{t,c} : \text{fin}(\mathcal{P}) \to \{0, 1\}$, where

$$H_{\mathcal{I}}^{t,c}(A) = \begin{cases} 1, & t(\{i \in \mathbb{N} : \mathcal{I}(i) = A\}) > c, \\ 0, & \text{otherwise,} \end{cases} \quad A \in \text{fin}(\mathcal{P}), \tag{16}$$

for some thresholding function $t$ and cutoff value $c \geq 0$. We define the *standard projection* $H_{\mathcal{I}}^*$ by taking $t(A) = \#A$ and $c = 0$ in (16).

We call $(\mathcal{Y}_n)_{n \geq 1}$ a *binary Hollywood process with parameter* $(\alpha, \theta)$ if it follows the Hollywood process with parameter $(\alpha, \theta, \nu)$ having $\nu_2 = 1$.

**Theorem 5.4** (Sparsity of projected network data). *Let $(\mathcal{Y}_n)_{n \geq 1}$ obey the binary Hollywood process with parameter $(\alpha, \theta)$. For each $n \geq 1$, let $H_n^{t,c}$ be the $(t, c)$-projection of $\mathcal{Y}_n$ by applying (16) for any $c \geq 0$ and $t(A) = \#A$, the cardinality map. Then the sequence $(H_n^{t,c})_{n \geq 1}$ is sparse almost surely for all $0 < \alpha < 1$ and $\theta > -\alpha$.*

Theorem 5.4 makes clear that projecting not only discards information but may also alter the assumed behavior of the data. By Theorems 5.3 and 5.4, $(\mathcal{Y}_n)_{n \geq 1}$ generated from the binary Hollywood model is sparse only for $1/2 < \alpha < 1$, while the projected network after applying (16) is sparse for all $0 < \alpha < 1$. We discuss these implications further in Section 6.4.

On the other hand, Figure 6 suggests that the power law behavior of the Hollywood model, as established in Theorem 5.2, might be preserved under projection. Whether this phenomenon is real or perceived, there remains no logical justification for projecting interaction data to a simple graph, especially when this operation makes the otherwise easy practice of parameter estimation intractable, as we discuss further surrounding (21) below. We also note that thresholding an edge exchangeable network according to (16) does not preserve edge exchangeability.

## 6. Inference from edge exchangeable models

6.1. **Maximum likelihood estimation.** Given edge-labeled network data $\mathcal{Y}_n$ with $n$ edges, the log-likelihood $\ell(\alpha, \theta, \nu; \mathcal{Y}_n)$ based on (13) with parameter $(\alpha, \theta, \nu)$ satisfies

$$\ell(\alpha, \theta, \nu; \mathcal{Y}_n) = \sum_{k=1}^{\infty} M_k(\mathcal{Y}_n) \log \nu_k + v(\mathcal{Y}_n) \log(\alpha) + \sum_{j=0}^{v(\mathcal{Y}_n)-1} \log(\theta/\alpha + j) - \\ - \sum_{j=1}^{m(\mathcal{Y}_n)} \log(\theta + j - 1) + \sum_{k=2}^{\infty} N_k(\mathcal{Y}_n) \sum_{j=0}^{k-2} \log(1 - \alpha + j), \tag{17}$$

where $M_k(\mathcal{Y}_n)$ is the number of edges in $\mathcal{Y}_n$ with exactly $k$ vertices, $N_k(\mathcal{Y}_n)$ is the number of vertices in $\mathcal{Y}_n$ with degree $k$, and $m(\mathcal{Y}_n) = \sum_{k=1}^{\infty} k M_k(\mathcal{Y}_n)$ is the total degree of $\mathcal{Y}_n$.

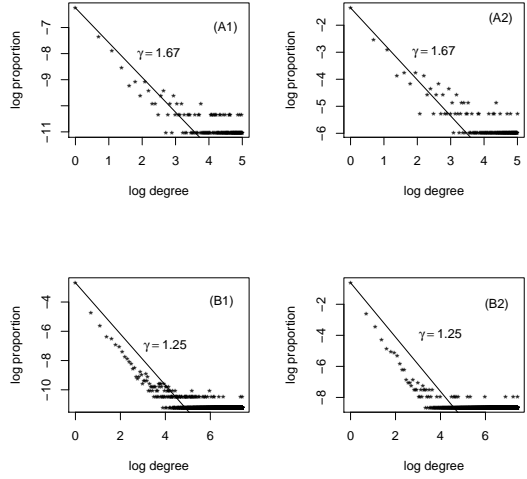

FIGURE 6. Simulated degree distribution of networks and their standard projection. Plots show the degree distribution of network obtained from (A1) the binary Hollywood model with $(\alpha, \theta) = (0.67, 1)$, (B1) the binary Hollywood model with $(\alpha, \theta) = (0.25, 1)$, (A2) the standard projection of the network in (A1), and (B2) the standard projection of the network in (B1). The line with slope $-\gamma$ in (A1) and (B1) indicates the true power law based on Theorem 5.2. The line with slope $-\gamma$ in (A2) and (B2) is the conjectured power law based on visual evidence and the connection to Theorem 5.2.

Maximum likelihood estimation of $\nu$ returns the empirical distribution $\hat{\nu}_{MLE} = \{\hat{\nu}_k\}_{k\geq 1}$, where $\hat{\nu}_k = M_k(\mathcal{Y}_n)/n$ for each $k = 1, 2, \ldots$. We then estimate $\alpha$ and $\theta$ by iterating between the score functions

$$(18) \quad \frac{\partial \ell(\alpha, \theta, \nu; \mathcal{Y}_n)}{\partial \alpha} = \frac{v(\mathcal{Y}_n)}{\alpha} + \sum_{j=0}^{v(\mathcal{Y}_n)-1} \frac{-\theta/\alpha^2}{\theta/\alpha + j} - \sum_{k=2}^{\infty} \sum_{j=0}^{k-2} \frac{N_k(\mathcal{Y}_n)}{1 - \alpha + j} = 0 \quad \text{and}$$

$$(19) \quad \frac{\partial \ell(\alpha, \theta, \nu; \mathcal{Y}_n)}{\partial \theta} = \sum_{j=0}^{v(\mathcal{Y}_n)-1} \frac{1/\alpha}{\theta/\alpha + j} - \sum_{j=0}^{m(\mathcal{Y}_n)-1} \frac{1}{\theta + j} = 0.$$

We encounter no convergence issues when iterating between (18) and (19) for our applications below.

The distribution in (13), and therefore the log-likelihood in (17), applies to the oriented network generated by the Hollywood model with parameter $(\alpha, \theta, \nu)$. With $\mathcal{Y}_n^*$ denoting the undirected edge-labeled network obtained by removing the orientations from $\mathcal{Y}_n$, for each $n \geq 1$, we have

$$\mathrm{pr}(\mathcal{Y}_n^* = \mathcal{E}^*; \alpha, \theta, \nu) =$$

$$(20) \quad C(\mathcal{E}^*) \left[ \prod_{k \geq 1} v_k^{M_k(\mathcal{E}^*)} \right] \alpha^{v(\mathcal{E}^*)} \frac{(\theta/\alpha)^{\uparrow v(\mathcal{E}^*)}}{\theta^{\uparrow m(\mathcal{E}^*)}} \prod_{k=2}^{\infty} \exp\{N_k(\mathcal{E}^*) \log((1-\alpha)^{\uparrow(k-1)})\},$$

where $C(\mathcal{E}^*)$ is the combinatorial factor that counts the number distinct ways to orient the edges of $\mathcal{E}^*$ to obtain a directed network. The log-likelihood based on $\mathcal{Y}_n^*$ satisfies

$$\tilde{\ell}(\alpha, \theta, \nu; \mathcal{Y}_n^*) = \ell(\alpha, \theta, \nu; \mathcal{Y}_n) + \log C(\mathcal{Y}_n^*),$$

where $\ell(\alpha, \theta, \nu; \mathcal{Y}_n)$ is the log-likelihood from (17) with $\mathcal{Y}_n$ taken to be any oriented edge-labeled network whose edges agree with $\mathcal{Y}_n^*$. Thus, the score functions based on $\mathcal{Y}_n^*$ are just as in (18) and (19) and maximum likelihood estimation for $(\alpha, \theta, \nu)$ can be performed without issue.

Problems arise, however, when projecting multiple edges according to the operation in (16). For definiteness, suppose $H_n^*$ is the standard projection from the binary Hollywood model with parameter $(\alpha, \theta)$ on $\mathfrak{E}_{[n]}$ as in (13). The log-likelihood $\ell^*(\alpha, \theta, \nu; H_n^*)$ in this case satisfies

$$\exp\{\ell^*(\alpha, \theta, \nu; H_n^*)\} =$$

$$= \sum_{\mathcal{E} \geq H_n^*} \left[ \prod_{k \geq 1} v_k^{M_k(\mathcal{E})} \right] \alpha^{v(\mathcal{E})} \frac{(\theta/\alpha)^{\uparrow v(\mathcal{E})}}{\theta^{\uparrow m(\mathcal{E})}} \prod_{k=2}^{\infty} \exp\{N_k(\mathcal{E}) \log((1-\alpha)^{\uparrow(k-1)})\}$$

$$= \left[ \prod_{k \geq 1} v_k^{M_k(H_n^*)} \right] \frac{\alpha^{v(H_n^*)} (\theta/\alpha)^{\uparrow v(H_n^*)}}{\theta^{\uparrow m(H_n^*)}} \prod_{k=2}^{\infty} \exp\{N_k(H_n^*) \log(1-\alpha)^{\uparrow(k-1)}\} \times$$

$$\times \sum_{\mathcal{E} \geq H_n^*} \frac{\prod_{k \geq 1} v_k^{M_k(\mathcal{E}) - M_k(H_n^*)}}{(\theta + m(H_n^*))^{\uparrow(m(\mathcal{E}) - m(H_n^*))}} \exp\{(N_k(\mathcal{E}) - N_k(H_n^*)) \log(1-\alpha)^{\uparrow(k-1)}\},$$

where $\mathcal{E} \geq H_n^*$ indicates that $\mathcal{E} \in \mathfrak{E}_{[n]}$ is an edge-labeled network whose standard projection is $H_n^*$. The log-likelihood has the form

$$(21) \quad \ell^*(\alpha, \theta, \nu; H_n^*) = \log C_n^*(H_n^*; \alpha, \theta, \nu) + \ell(\alpha, \theta, \nu; H_n^*),$$

so that without a manageable expression for $C_n^*(H_n^*; \alpha, \theta, \nu)$, maximum likelihood estimation based on projected network data $H_n^*$ is intractable.

6.2. **Application to Wikipedia voting and Karate Club networks.** Table 1 shows the maximum likelihood estimates for the binary Hollywood model fit to the Wikipedia voting and Karate Club networks described in Crane and Dempsey (2016, Section B). Together these examples cover both regimes in the parameter space of the Hollywood model.

The Karate Club network consists of interactions among thirty-four club members, warranting the choice of $k = 34$ when fitting the model with parameters $\alpha < 0$ and $\theta = -34\alpha$. The Wikipedia network, by contrast, has

Maximum likelihood estimates

| network | $\hat{\alpha}_{MLE}$ | std. error* | $\hat{\theta}_{MLE}$ | std. error* |
|---|---|---|---|---|
| Wikipedia | 0.37 | 0.01 | 183 | 14.76 |
| Karate Club | -1.80 | 0.47 | 61.3 | 16.04 |

TABLE 1. Maximum likelihood estimates and standard errors for $(\alpha, \theta)$ from binary Hollywood model fit to Wikipedia voting and Karate Club networks. * Standard errors estimated by the Cramér–Rao lower bound, which we have verified as a good estimate for the standard error based on simulation.

no upper limit on the number of vertices, and so we fit the model under regime $0 < \alpha < 1$ and $\theta > -\alpha$.

The large standard error for maximum likelihood estimates of $\theta$ agrees with what is known about estimation of the mutation rate in Ewens's sampling formula (Crane, 2016a; Ewens, 1972): although $\hat{\theta}_{MLE} \to \theta$ almost surely as the sample size grows, it converges at a rate on the order of $\log(n)$, rendering it inconsistent for practical purposes. The estimate of $\alpha$ is of greater interest in the applications we envision because of its relationship to the power law behavior, cf. Theorem 5.2. When $\alpha < 0$ and the total number of vertices is taken to be known and finite $k = 1, 2, \ldots$, then $\theta = -k\alpha$ and we obtain $\hat{\theta}_{MLE} = k\hat{\alpha}_{MLE}$. In this case, we have s.e.$(\hat{\theta}_{MLE}) = k \times$ s.e.$(\hat{\alpha}_{MLE})$, with s.e.$(\cdot)$ shorthand for the standard error.

By Theorem 5.2, the maximum likelihood estimate $\hat{\alpha}_{MLE} = 0.37$ for the Wikipedia dataset implies an estimated power law of $\hat{\gamma}_{MLE} = 1.37$, which is reasonably close to the estimated power law exponent of $\hat{\gamma}_{YULE} = 1.44$ obtained by fitting degree distribution of the network to the Yule model,

(22) $$\text{pr}(K = k; \gamma) = (\gamma - 1)\Gamma(\gamma)/\Gamma(k + \gamma), \quad k \geq 1,$$

for $\gamma > 1$. The likelihood surface is continuous in $\alpha, \theta$, except at $\alpha = 0$, and estimation experienced no instability or convergence issues in multiple iterations from different starting points.

6.3. **Prediction using growth dynamics.** The growth dynamics of edge exchangeable models by the sequential process of edge addition facilitates predictive inferences in networks generated by a process of repeated interactions. In the setting of Section 5, we can predict the next interaction, based on an observed network $\mathcal{Y}_n$, from the update probabilities in (12) with $(\alpha, \theta, \nu) = (\hat{\alpha}_{MLE}, \hat{\theta}_{MLE}, \hat{\nu}_{MLE})$ given by the maximum likelihood estimates obtained from (17).

For a concrete application, we consider the actors collaboration network from Barabási and Albert (1999). For $(\alpha, \theta, \nu)$ in the parameter space of the Hollywood model, we compute the predictive probability, conditional on $\mathcal{Y}_n$, that the next movie contains an actor not previously observed in the sample. A straightforward calculation based on (12) and the law of total

Maximum likelihood estimates

| k | $\hat{v}_k$ | k | $\hat{v}_k$ |
|---|---|---|---|
| 1 | 0.081 | 7 | 0.057 |
| 2 | 0.071 | 8 | 0.054 |
| 3 | 0.065 | 9 | 0.046 |
| 4 | 0.072 | 10 | 0.275 |
| 5 | 0.062 | 11 | 0.161 |
| 6 | 0.059 | | |

TABLE 2. Maximum likelihood estimates of the movie size distribution for the actors collaboration network. The mean number of actors in each movie is 7.136.

probability gives

$$(23) \quad \text{pr}(\text{new actor in next movie} \mid \mathcal{Y}_n; \alpha, \theta, v) = 1 - \sum_{k \geq 1} v_k \frac{(M - \alpha N)^{\uparrow k}}{(\theta + M)^{\uparrow k}},$$

where $N = \sum_{k \geq 1} N_k(\mathcal{Y}_n)$ is the total degree of $\mathcal{Y}_n$ and $M = \sum_{k \geq 1} k N_k(\mathcal{Y}_n)$ is the total number of roles in $\mathcal{Y}_n$. (To see the calculation in (23), we note that $(M - \alpha N)^{\uparrow k}/(\theta + M)^{\uparrow k}$ is the conditional probability that the next movie does not feature a new actor, given that it has $k$ actors in its cast. Summing over all $k \geq 1$ gives the total probability that the next movie does feature a new actor, yielding the rightmost term in (23). The probability that a new actor appears follows by taking the complementary probability.)

Fitting the model to the actors collaboration dataset yields $\hat{\alpha}_{MLE} = 0.66$ (s.e. $6.8 \times 10^{-4}$) and $\hat{\theta}_{MLE} = 4.21$ (s.e. 2.86), with $\hat{v}_{MLE}$ given in Table 2. Since $\hat{\alpha}_{MLE} = 0.66$ lies in the range $1/\mu < \alpha < 1$, with $\mu = 7.136$ computed from Table 2, Theorem 5.3 suggests that the actors collaboration network is sparse. This estimate also agrees exactly with the estimated fit $\hat{\gamma}_{YULE} = 1.66$ if the degree distribution is fit directly to the Yule distribution in (22).

The estimated predictive probability based on these maximum likelihood estimates is 0.78. We check the accuracy of this estimate by data-splitting cross validation based on samples of 2,000 movies from the total collection of about 200,000. For each iteration $j = 1, 2, \ldots$, we sample 2,000 movies uniformly at random from the 200,000 to obtain $\mathcal{Y}_{2000}$. We then compute the probability in (23) based on the estimates $\hat{\alpha}_{MLE}$, $\hat{\theta}_{MLE}$, and $\hat{v}_{MLE}$ from $\mathcal{Y}_{2000}$ and compare this estimated probability to the empirical probability computed as the proportion of the 198,000 unsampled movies for which there appears an actor not among the 2,000 movies in the sample $\mathcal{Y}_{2000}$. The mean relative error between the fitted probability (23) and the empirical probability for 100 iterations was $-0.003$ with a standard deviation of 0.002.

6.4. **Tests for sparsity and power law.** As asymptotic network properties, sparsity and power law cannot be verified with certainty based on any finite amount of data. Statistical tests for sparsity and/or power law based on finite sample data, therefore, require that the finite sample models faithfully represent the properties exhibited by the population network.

Our discussion in Section 5.4 advises caution when testing for asymptotic properties based on thresholded network data. If the parameter space $\Pi$ partitions as $\Pi = \Pi_0 \cup \Pi_1$ such that $\pi \in \Pi_1$ parameterizes sparse networks and $\pi \in \Pi_0$ parameterizes non-sparse networks, then $H_0 : \pi \in \Pi_0$ versus $H_1 : \pi \in \Pi_1$ might yield a valid test for sparsity. Alternative interpretations of the partition $\Pi = \Pi_0 \cup \Pi_1$, however, may provide a more parsimonious conclusion. For example, under the binary Hollywood model of Section 5, the region $1/2 < \alpha < 1$ corresponds to sparse networks while $\alpha < 1/2$ parameterizes dense networks, with $\alpha < 0$ corresponding to the case where the number of vertices stays bounded and finite. By Theorem 5.4, if the network is thresholded by the standard projection, then the projected network is sparse a.s. for all $\alpha > 0$. In this case, the partition $\Pi = \{\alpha < 0\} \cup \{0 \leq \alpha \leq 1\}$ offers the more parsimonious interpretation according to whether the population of vertices is finite ($\alpha < 0$) or infinite ($0 \leq \alpha \leq 1$) rather than sparse ($0 \leq \alpha \leq 1$) or not sparse ($\alpha < 0$).

A particular consequence of this is on display when testing for sparsity in the US airport dataset (Opsahl, 2011), which is built from the flight map between all US airports in 2010. For each pair of airports, there are as many edges as there were seats on all flights between those airports. The edge weights range over several orders of magnitude from a minimum of 1 to several hundred thousand, making the projection in (16) particularly deleterious to the data structure. For the binary Hollywood model fit to this dataset, there are two possible maximum likelihood estimates, both with log-likelihood of about $-8.21 \times 10^{-9}$: if $\alpha > 0$ we get $\hat{\alpha}_{MLE,1} = 0.13$ (s.e. 0.003) and $\hat{\theta}_{MLE,1} = 0.08$ (s.e. 0.134), and if $\alpha < 0$ we get $\hat{\alpha}_{MLE,2} = -0.11$ (s.e. 0.003) by taking $k = 1574$ the number of airports in the sample. The estimate $\hat{\alpha}_{MLE,1} = 0.13$ is consistent with $\hat{\gamma}_{YULE} = 1.11$ obtained by fitting (22) to the degree distribution of the airport network. The choice of regime $\alpha < 0$ or $\alpha > 0$ is a matter of whether we believe the population of airports under consideration is fixed and finite or potentially unbounded, a choice which is closely tied to the desired hypothesis test.

Technically, it makes sense to test for sparsity only if the population is infinite, in which case $\alpha > 0$ and we can test $H_0 : 0 < \alpha \leq 1/2$ against $H_1 : 1/2 < \alpha < 1$. From the estimate $\hat{\alpha}_{MLE,1} = 0.13$, we cannot reject the hypothesis $H_0 : \alpha < 1/2$. But if testing based on the projected network, then the above test is degenerate since the ranges for which the population is infinite and for which the population is sparse coincide, per Theorem 5.4. By contrast, Caron and Fox (2014, Section 7.2) conclude that the US airport network is sparse based on a 99% credible interval of $[0.099, 0.181]$ for an analogous parameter to $\alpha$ in our case, but that analysis appears to be based on the standard projection.

## 7. Discussion of other approaches

The above discussion offers ample reason to prefer the edge exchangeable approach over other common models for network datasets constructed from interaction processes. The observation that the interactions, and therefore edges, comprise the units leads naturally to the alternative representation by edge-labeled networks. Beyond treating the units appropriately, representing the data by an edge-labeled network frees us of the limitations of the conventional representation by vertex-labeled graphs. The Hollywood model, in particular, is computationally tractable, performs admirably in several real data examples, and replicates the key features of sparsity and power law degree distribution of primary interest in the modern network science literature. The nonparametric vertex components model in Section 4.2 offers an even larger class of models with likely computational benefits due to its stick-breaking construction.

To round out the discussion, we compare the performance of the Hollywood model to that of the exponential random graph model (ERGM) (Holland and Leinhardt, 1981) and stochastic blockmodel (SBM) (Holland et al., 1983; Snijders and Nowicki, 1997) fit to the actors collaboration network. We then discuss some treatments of sparse networks by other authors.

7.1. **Empirical comparison.** We compare the performance of the Hollywood model (HW), stochastic blockmodel (SBM), Erdős–Rényi model (ER), and exponential random graph model (ERGM) on the actors collaboration dataset (Barabási and Albert, 1999), which consists of approximately $200,000$ movies. We have already fit HW to the full movie dataset in Section 6.3, but we could not compare to other methods on the full dataset because the best available software for SBM[2] and ERGM[3] are computationally inefficient for networks of even moderate size.

For comparison, we subsample 250 movies uniformly at random from the database of $200,000$ movies. Since SBM, ER, and ERGM are models for graphs without multiple edges, we treat the *de facto* population network as the graph $G_{250}$ induced by these 250 movies as in (1). The maximum likelihood estimate (MLE) under SBM estimates 21 blocks and log-likelihood of $-32,818$ and MLE under ER returns a log-likelihood of $-42,924$. The ERGM with sufficient statistics given by edge and triangle density, respectively,

$$\varepsilon(G) = \frac{2}{v(G)(v(G)-1)} e(G) \quad \text{and} \quad \tau(G) = \frac{6}{v(G)(v(G)-1)(v(G)-2)} \Delta(G),$$

where $\Delta(G)$ is the number of triangles in $G$, returns a model degeneracy error.

We also fit ERGM and ER to the induced graph (1) obtained from a $1,000$ movie subsample. In this case, the fitted parameters for ERGM are $\hat{\beta}_\varepsilon = -12.26$ (s.e. 0.083) and $\hat{\beta}_\tau = 5.24$ (s.e. 0.019) with log-likelihood of $-190,630$, which improves over Erdős–Rényi (log-likelihood: $-208,702$). By contrast, we fit the Hollywood model to both the 250 and 1,000 movie subsamples and obtain a maximum log-likelihood of $-476$ and $-26,320$, respectively.

---

[2] https://cran.r-project.org/web/packages/blockmodels/blockmodels.pdf
[3] https://cran.r-project.org/web/packages/ergm/ergm.pdf

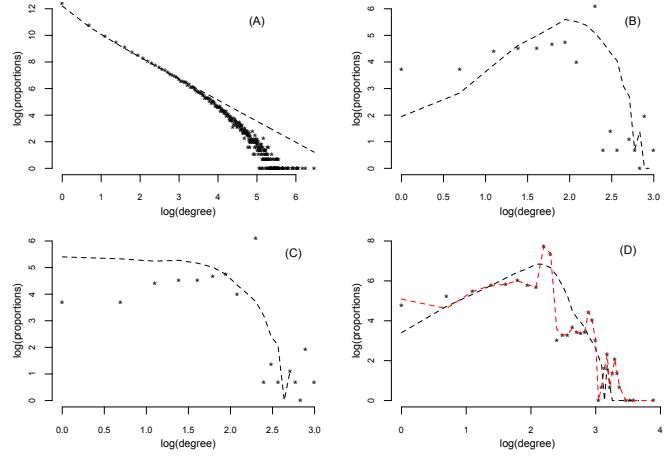

FIGURE 7. Comparison of fitted degree distribution (- - -) to empirical degree distribution (*) for (A) Hollywood model (compared to full dataset), (B) Erdős–Rényi model and (C) stochastic blockmodel (compared to 250 movie sample), and (D) ERGM (red) and ER (black) (compared to 1,000 movie sample).

Figure 7 compares the goodness of fit for the degree distribution for HW, SBM, ER, and ERGM. The fitted degree distribution for HW in (A) is based on the estimate obtained from the sample of 250 movies compared to the degree distribution of the complete population network on 200,000. In (B) and (C), we compare the fitted ER and SBM to the degree distribution of subsampled network $G_{250}$ induced by (1) and in (D) we compare the fitted ER and ERGM to the degree distribution of $G_{1000}$ based on 1,000 movies through (1). In this sense, the comparison is generous to SBM, ER, and ERGM since those are compared only to the data to which they were fit while HW is evaluated based on how well it detects the population structure for 200,000 edges based on a sample of 250 movies. Even with this disadvantage, the Hollywood model gives a much better fit than ER and SBM, while the fit under ERGM is very close to the observed degree distribution and shows strong evidence of overfitting. Finally, the Hollywood model converged to the MLE without problem while SBM and ERGM experienced convergence issues and long run times.

### 7.2. Further discussion.

7.2.1. *Preferential attachment models.* The most well known model for sparse, power law networks is the *preferential attachment model* of Barabási and Albert (1999), which successfully puts forward a generating mechanism for sparse networks that evolve by vertex growth and have a power law degree distribution of exponent $\gamma > 2$. The model, however, seems suitable as a statistical model only for specially structured datasets, such as citation networks (Simon, 1955; De Solla Price, 1965). Otherwise, if the vertex ordering is unobserved or there is no natural ordering of the vertices, then these models lack basic logical properties, namely *label equivariance*, to ensure sound inferences. Kolaczyk (2009, Section 6.4.3) also discusses some practical issues with fitting the preferential attachment model to data.

7.2.2. *Graphon models.* Graphon models grow out of the theory of dense graph limits (Lovász and Szegedy, 2006) and are intimately related to the Aldous–Hoover theory of partially exchangeable arrays (Aldous, 1981; Hoover, 1979). As we mentioned in Section 1, traditional exchangeable random graph models with vertex set $\mathbb{N}$ cannot reproduce sparsity.

Probabilists often address this issue by studying the graphon under a *sparse regime*, an approach brought into the statistics literature by Bickel and Chen (2009). With $\phi : [0,1] \times [0,1] \to [0,1]$ symmetric in its arguments and $\{\rho_n\}_{n \geq 1}$ satisfying $\rho_n \to \infty$ and $\rho_n^{-1} n \int_{[0,1]\times[0,1]} \phi(u,v) du\, dv = O(1)$ as $n \to \infty$, Bickel and Chen (2009) define the distribution of a graph $G_n = ([n], E_n)$ with vertices labeled in $[n]$ by taking $U_1, \ldots, U_n$ i.i.d. Uniform$[0,1]$ and putting

(24) $\quad \text{pr}(\{i,j\} \in E_n \mid U_1, \ldots, U_n) = \rho_n^{-1} \phi(U_i, U_j), \quad 1 \leq i < j \leq n,$

conditionally independently for all $1 \leq i < j \leq n$. Under (24), every $G_n$ is exchangeable, and collectively the sequence $(G_n)_{n=1,2,\ldots}$ is sparse in the sense of Definition 3.1. But the formulation does not automatically correspond to a generating process for a population network since the marginal distributions in (24) are not logically related through any specified sampling mechanism.

Without a direct link between samples of different size, the meaning of the parameter $\phi$ varies with sample size and there is no logical way to relate parameter estimates $\hat{\phi}_n$ based on a sample of size $n$ to a single parameter $\phi$ for the population. In this case, the interpretation of estimated parameters is obscured and the significance of asymptotic statements for finite sample approximations, as derived, for example, by Borgs et al. (2015) and Wolfe and Olhede (2014), is unclear. Shalizi and Rinaldo (2013) discuss related issues for estimation from exponential random graph models. Crane and Dempsey (2015) discuss the effects of sampling more generally.

7.2.3. *Completely random measures and graphex models.* Caron and Fox (2014), and later Veitch and Roy (2015) and Borgs et al. (2016), studied a class of network models derived from exchangeable point processes $\mathbf{X}$ on $[0,\infty)\times[0,\infty)$. In this setting, each $(x, x') \in \mathbf{X}$ determines an edge between vertices corresponding to points $x, x'$. The point process $\mathbf{X}$ is assumed to be *exchangeable* in the sense that its distribution is invariant under the joint action on $[0, \infty) \times [0, \infty)$ by measure preserving transformations of $[0, \infty)$.

In this context, a sequence of network data is obtained by defining a graph $G_t$ for each $t \geq 0$, where $G_t$ is derived from the restriction of $\mathbf{X}$ to $[0, t] \times [0, t]$ by only including in $G_t$ those vertices that are labeled in $[0, t]$ and which are not isolated in the restriction. Veitch and Roy (2015) adapt a theorem of Kallenberg (1990) to obtain a generic *graphex representation* of the class of random graphs derived from exchangeable point processes in this way.

Graphex models allow for sparse, power law behavior, and their description in terms of the exchangeable point process $\mathbf{X}$ does invoke some invariance principle in the construction. It is important to stress, however, that the notion of exchangeability applies to the generating point process $\mathbf{X}$, and not directly to the associated process of graphs $(G_t)_{t \geq 0}$. In fact, the graphs $(G_t)_{t \geq 0}$ induced from the process are not equipped with any labels, and so it is not yet clear what exchangeability of the point process implies for the induced graph sequence, or how the construction of $G_t$ by restricting to $[0, t]$ and removing isolated vertices relates to the manner in which real world networks are observed. At the time of this writing, graphex models are still under development with some of these questions in mind.

## 8. Proof of main theorems

We equip $\mathfrak{E}_\mathbb{N}$ with the topology and Borel $\sigma$-field induced by the metric

$$d(\mathcal{E}, \mathcal{E}') = 1/(1 + \sup\{n \in \mathbb{N} : \mathcal{E}|_{[n]} = \mathcal{E}'|_{[n]}\}), \quad \mathcal{E}, \mathcal{E}' \in \mathfrak{E}_\mathbb{N},$$

with the convention $1/\infty = 0$.

8.1. **Proof of Theorem 5.1.** Let $(\mathcal{Y}_n)_{n \in \mathbb{N}}$ be a realization of the Hollywood process with parameter $(\alpha, \theta, \nu)$ from Section 5. Then the distribution of $\mathcal{Y}_n$ is given by (13) for each $n \in \mathbb{N}$, and these distributions are consistent with respect to the restriction operation by consequence of the sequential construction in (12). The distribution in (13) depends on $\mathcal{Y}_n$ only through $v(\mathcal{Y}_n)$, $N_k(\mathcal{Y}_n)$, and $M_k(\mathcal{Y}_n)$, that is, the number of vertices, number of vertices of degree $k \geq 1$, and the number of $k$-ary edges for each $k \geq 1$, all of which are invariant under relabeling edges. Edge exchangeability follows.

8.2. **Proof of Theorem 5.2.** Let $(\mathcal{Y}_n)_{n \geq 1}$ be a realization from the Hollywood process with parameter $(\alpha, \theta, \nu)$. The power law behavior is apparent by the following connection between $(\mathcal{Y}_n)_{n \geq 1}$ and the two parameter Ewens–Pitman distribution. In particular, the Ewens–Pitman distribution corresponds to the special case of the Hollywood model with $0 < \alpha < 1$, $\theta > -\alpha$, and $\nu_1 = 1$. In this case, there is exactly one vertex incident to every edge and the network data corresponds to a partition of the edges. The number of blocks in the Ewens–Pitman partition corresponds to the number of vertices in the unary Hollywood process. Let $N_k(\mathcal{Y}_n)$ be the number of vertices with degree $k$ in $\mathcal{Y}_n$ and let $v(\mathcal{Y}_n)$ be the total number of vertices in $\mathcal{Y}_n$. By Pitman (2005, Lemma 3.11), we know that $N_k(\mathcal{Y}_n)/v(\mathcal{Y}_n) \to \alpha(1-\alpha)^{\uparrow(k-1)}/k!$

a.s. for every $k \geq 1$ as $n \to \infty$. Thus, the sequence of degree distributions $\{d_n\}_{n\geq 1} = \{(N_k(\mathcal{Y}_n)/v(\mathcal{Y}_n))_{k\geq 1}\}_{n\geq 1}$ converges a.s. in the total variation topology to the distribution given by $p_\alpha(k) = \alpha(1-\alpha)^{\uparrow(k-1)}/k!$.

In the general Hollywood process with arbitrary distribution $\nu$, the degree sequences of $(\mathcal{Y}_n)_{n\geq 1}$ coincide with a random subsequence of $\{d_n\}_{n\geq 1}$ indexed by $\{K_r\}_{r\geq 1}$ for $K_r = \sum_{j=1}^r \kappa_j$, $r = 1, 2, \ldots$, where $\kappa_1, \kappa_2, \ldots$ are i.i.d. from $\nu$. Thus, the degree distributions of $(\mathcal{Y}_n)_{n\geq 1}$ correspond to $\{d_{K_r}\}_{r\geq 1}$, which is a subsequence of the a.s. converging sequence $\{d_n\}_{n\geq 1}$ and, therefore, must have the same a.s. limit. The proof is completed by noting that $\alpha(1-\alpha)^{\uparrow(k-1)}/k! \sim \alpha k^{-(\alpha+1)}/\Gamma(1-\alpha)$ as $k \to \infty$, which corresponds to the power law with exponent $\alpha + 1$.

8.3. **Proof of Theorem 5.3.** We once again exploit the connection to the Ewens–Pitman distribution from the proof of Theorem 5.2. Let $N_n$ be the number of vertices in the unary Hollywood process with parameter $(\alpha, \theta)$ satisfying $0 < \alpha < 1$ and $\theta > -\alpha$. Pitman (2005, Theorem 3.8) shows that $n^{-\alpha}N_n \to S_\alpha$ a.s., where $S_\alpha$ is a strictly positive and finite random variable. The sequence $(N_n)_{n\geq 1}$ then satisfies $N_n \sim n^\alpha S_\alpha$ a.s. as $n \to \infty$, so that

$$E(N_n) \sim \frac{\Gamma(\theta+1)}{\alpha\Gamma(\theta+\alpha)} n^\alpha \quad \text{as } n \to \infty,$$

by Pitman (2005, Theorem 3.8).

Let $(\mathcal{Y}_n)_{n\geq 1}$ be the Hollywood process with parameter $(\alpha, \theta, \nu)$ for arbitrary distribution $\nu$. Then $(v(\mathcal{Y}_n))_{n\geq 1}$ is a random subsequence of $(N_n)_{n\geq 1}$ given by $(N_{K_r})_{r\geq 1}$, where $K_r = \sum_{j=1}^r \kappa_j$ for $\kappa_1, \kappa_2, \ldots$ i.i.d. from $\nu$. It follows from the above argument that $v(\mathcal{Y}_n) = N_{K_n} \sim K_n^\alpha S_\alpha$ a.s. as $n \to \infty$. By the strong law of large numbers, $n^{-1}K_n \sim \mu$ a.s. and $K_n^\alpha \sim (\mu n)^\alpha$ a.s. as $n \to \infty$, where $\mu = \sum_{k\geq 1} k\nu_k$; whence,

$$E(v(\mathcal{Y}_n)) \sim \frac{\Gamma(\theta+1)}{\alpha\Gamma(\theta+\alpha)} (\mu n)^\alpha \quad \text{as } n \to \infty.$$

To establish sparsity of $(\mathcal{Y}_n)_{n\geq 1}$, we consider $(n^{-1}v(\mathcal{Y}_n)^{m_\bullet(\mathcal{Y}_n)})_{n\geq 1}$, where $m_\bullet(\mathcal{Y}_n)$ is the average total degree in $\mathcal{Y}_n$. We must identify values of $0 < \alpha < 1$ for which $\liminf_{n\to\infty} n^{-1}v(\mathcal{Y}_n)^{m_\bullet(\mathcal{Y}_n)} = +\infty$. By the strong law of large numbers, $m_\bullet(\mathcal{Y}_n) \to \mu$ a.s. as $n \to \infty$. By the above discussion, $v(\mathcal{Y}_n) \sim (\mu n)^\alpha S_\alpha$ a.s. as $n \to \infty$ for a strictly positive and finite random variable $S_\alpha$. It follows that $n^{-1}v(\mathcal{Y}_n)^{m_\bullet(\mathcal{Y}_n)} \sim n^{-1}(\mu n)^{\mu\alpha} S_\alpha$ a.s. as $n \to \infty$, which goes to infinity as long as $\mu\alpha > 1$. Thus, $(\mathcal{Y}_n)_{n\geq 1}$ is sparse with probability 1 provided $1/\mu < \alpha < 1$, for all $\theta > -\alpha$.

8.4. **Proof of Theorem 5.4.** Let $(\mathcal{Y}_n)_{n\geq 1}$ be the binary Hollywood process with parameter $0 < \alpha < 1$ and $\theta > -\alpha$. For each $n \geq 1$, let $H_n^{t,c}$ be the projection obtained by applying (16) to $\mathcal{Y}_n$. To establish sparsity of $(H_n^{t,c})_{n\geq 1}$,

we must show
$$\limsup_{n\to\infty} \frac{e(H_n^{t,c})}{v(H_n^{t,c})^2} = 0 \quad \text{a.s.,}$$
for which it suffices to establish sparsity under the standard projection $H_n^*$, since $e(H_n^{t,c}) \leq e(H_n^*)$ and $v(H_n^{t,c}) = v(H_n^*)$.

For each $k \geq 1$ and $n \geq 1$, we write $N_k(\mathcal{Y}_n)$ and $N_k(H_n^*)$ to denote the number of vertices with degree $k$ in $\mathcal{Y}_n$ and $H_n^*$, respectively. We also write $e(\mathcal{Y}_n)$ and $e(H_n^*)$ to denote the number of edges in $\mathcal{Y}_n$ and $H_n^*$, respectively, and $v(\mathcal{Y}_n) = v(H_n^*)$ to denote the number of vertices in $\mathcal{Y}_n$ and $H_n^*$. Since the projection operation reduces multiple occurrences of the same edge to a single edge, the degree of each vertex in $H_n^*$ can be no larger than $v(H_n^*)$ and

$$(25) \quad e(H_n^*) = \sum_{k\geq 1} k N_k(H_n^*) \leq \sum_{k\geq 1} (k \wedge v(H_n^*)) N_k(\mathcal{Y}_n).$$

For every $K \geq 1$, (25) implies

$$\limsup_{n\to\infty} \frac{e(H_n^*)}{v(H_n^*)^2} \leq \limsup_{n\to\infty} \sum_{k=1}^{K} \frac{k \wedge v(H_n^*)}{v(H_n^*)} \frac{N_k(\mathcal{Y}_n)}{v(H_n^*)} + \limsup_{n\to\infty} \sum_{k=K+1}^{\infty} \frac{N_k(\mathcal{Y}_n)}{v(H_n^*)}$$

$$\leq \limsup_{n\to\infty} \sum_{k=1}^{K} \frac{k \wedge v(H_n^*)}{v(H_n^*)} + \limsup_{n\to\infty} \sum_{k=K+1}^{\infty} \frac{N_k(\mathcal{Y}_n)}{v(H_n^*)}$$

$$(26) \quad \leq \sum_{k=1}^{K} \limsup_{n\to\infty} \frac{k \wedge v(H_n^*)}{v(H_n^*)} + \limsup_{n\to\infty} \sum_{k=K+1}^{\infty} \frac{N_k(\mathcal{Y}_n)}{v(H_n^*)}.$$

Pitman (2005, Corollary 3.9) implies that $n^{-\alpha} v(H_n^*) \to S_\alpha$ a.s., where $S_\alpha$ is a strictly positive, finite random variable; thus, $v(H_n^*) \to \infty$ a.s. and $\limsup_{n\to\infty}(k \wedge v(H_n^*))/v(H_n^*) \to 0$ a.s. for every $k = 1,\ldots,K$, so that (26) implies

$$(27) \quad \limsup_{n\to\infty} \frac{e(H_n^*)}{v(H_n^*)^2} \leq \limsup_{n\to\infty} \sum_{k=K+1}^{\infty} \frac{N_k(\mathcal{Y}_n)}{v(H_n^*)} \quad \text{for all } K \geq 1.$$

By Pitman (2005, Lemma 3.11),
$$\lim_{n\to\infty} \frac{N_k(\mathcal{Y}_n)}{v(H_n^*)} = p_\alpha(k) = \alpha \Gamma(k-\alpha)/(k!\Gamma(1-\alpha)) \quad \text{for every } k \geq 1 \text{ a.s.,}$$
implying that for every $\varepsilon, \delta > 0$ and $k \geq 1$ there exists $R = R(\varepsilon, \delta, k)$ such that
$$\mathrm{pr}(|v(H_n^*)^{-1} N_k(\mathcal{Y}_n) - p_\alpha(k)| < \varepsilon \text{ for all } n > R) \geq 1 - \delta.$$
For any $K \geq 1$, we choose $R^* = \max_{1 \leq k \leq K} R(\varepsilon/K, \delta/K, k)$ so that
$$\mathrm{pr}(|v(H_n^*)^{-1} N_k(\mathcal{Y}_n) - p_\alpha(k)| < \varepsilon/K \text{ for all } n > R^*, \text{ for all } k = 1,\ldots,K) \geq 1 - \delta$$
and
$$\mathrm{pr}\left(\left|\sum_{k=K+1}^{\infty} v(H_n^*)^{-1} N_k(\mathcal{Y}_n) - \sum_{k=K+1}^{\infty} p_\alpha(k)\right| < \varepsilon \text{ for all } n > R^*\right) \geq 1 - \delta.$$

We combine this with (27) and the tail calculation $p_\alpha(> K) = \sum_{k=K+1}^{\infty} p_\alpha(k) = \Gamma(K + 1 - \alpha)/(\Gamma(K+1)\Gamma(1-\alpha))$ to observe

$$\text{pr}(\limsup_{n\to\infty} e(H_n^*)/v(H_n^*)^2 > p_\alpha(> K) + \varepsilon) \leq$$

$$\leq \text{pr}\left(\limsup_{n\to\infty} \sum_{k=K+1}^{\infty} N_k(\mathcal{Y}_n)/v(H_n^*) > p_\alpha(> K) + \varepsilon\right)$$

$$\leq \text{pr}\left(\left|\sum_{k=K+1}^{\infty} v(H_n^*)^{-1} N_k(\mathcal{Y}_n) - \sum_{k=K+1}^{\infty} p_\alpha(k)\right| > \varepsilon \text{ for some } n > R^*\right)$$

$$\leq \delta$$

for every $K \geq 1$. Since $p_\alpha(> K) \downarrow 0$ as $K \to \infty$, it follows that

$$\text{pr}(\limsup_{n\to\infty} e(G_n^*)/v(H_n^*)^2 > \varepsilon) \leq \delta \quad \text{for all } \varepsilon, \delta > 0$$

and $(H_n^*)_{n\geq 1}$ is sparse a.s.

DEPARTMENT OF STATISTICS & BIOSTATISTICS, RUTGERS UNIVERSITY, 110 FRELINGHUYSEN AVENUE, PISCATAWAY, NJ 08854, USA
 *E-mail address*: hcrane@stat.rutgers.edu
 *URL*: http://stat.rutgers.edu/home/hcrane

DEPARTMENT OF STATISTICS, UNIVERSITY OF MICHIGAN, 1085 S. UNIVERSITY AVE, ANN ARBOR, MI 48109, USA
 *E-mail address*: wdem@umich.edu

# SUPPLEMENTARY MATERIAL: EDGE EXCHANGEABLE MODELS FOR NETWORK DATA

HARRY CRANE AND WALTER DEMPSEY

ABSTRACT. This is the supplement to Crane and Dempsey (2016). We characterize edge exchangeable network models and describe additional examples of interaction datasets.

## Appendix A. Characterization of edge exchangeable networks

Our representation in Crane and Dempsey (2016, Theorem 4.2) applies only to the special case of edge exchangeable networks that are *blip-free*. We prove here a more general, complete characterization for all edge exchangeable networks in the binary case for which each edge is undirected and involves exactly two vertices. The argument for directed networks with edges of any finite arity is analogous but more technical. Recall the notation for edge-labeled networks introduced in Crane and Dempsey (2016).

Let $\text{fin}_2(\mathcal{P}) \subset \text{fin}(\mathcal{P})$ be the set of all size 2 multisets of $\mathcal{P}$. Given a binary edge-labeled network $\mathcal{E}$ with $e(\mathcal{E}) = n$, we call $S : [n] \to \text{fin}_2(\mathbb{N})$ a *selection function* for $\mathcal{E}$ if $\mathcal{E}_S = \mathcal{E}$, where $\mathcal{E}_S$ is as defined by Crane and Dempsey (2016, Equation (5)); that is, $S$ is an interaction process whose induced edge-labeled network agrees with $\mathcal{E}$. We think of the selection function as a way of labeling the vertices of the edge-labeled network $\mathcal{E}$.

Selection functions $S, S' : [n] \to \text{fin}_2(\mathbb{N})$ are equivalent, written $S \equiv S'$, if they correspond to the same edge-labeled network, that is, $\mathcal{E}_S = \mathcal{E}_{S'}$. To every edge-labeled network $\mathcal{E}$ with $n$ edges we associate a *canonical selection function* $S_\mathcal{E} : [n] \to \text{fin}_2(\mathbb{N})$ defined by labeling the vertices *in order of appearance*, as follows.

We initialize by putting $S_\mathcal{E}(1) = \{1,1\}$ if the edge labeled 1 is a self loop and otherwise we put $S_\mathcal{E}(1) = \{1,2\}$ to indicate an edge between two distinct vertices. Given $S_\mathcal{E}(1), \ldots, S_\mathcal{E}(i-1)$, we define $S_\mathcal{E}(i) = \{v_1(i), v_2(i)\}$ by choosing $v_1(i) \le v_2(i)$ to be the smallest vertex labels consistent with the structure of $\mathcal{E}|_{[i]}$. In this way, $v_j(i) = r$ coincides with a previously observed vertex label if one of the vertices involved in the $i$th interaction corresponds to the vertex labeled $r$ in the previous interactions $S_\mathcal{E}(1), \ldots, S_\mathcal{E}(i-1)$. See Figures 1(b) and 1(c) for an illustration.

---

*Date*: October 19, 2016.
H. Crane is partially supported by NSF grants CNS-1523785 and CAREER DMS-1554092.

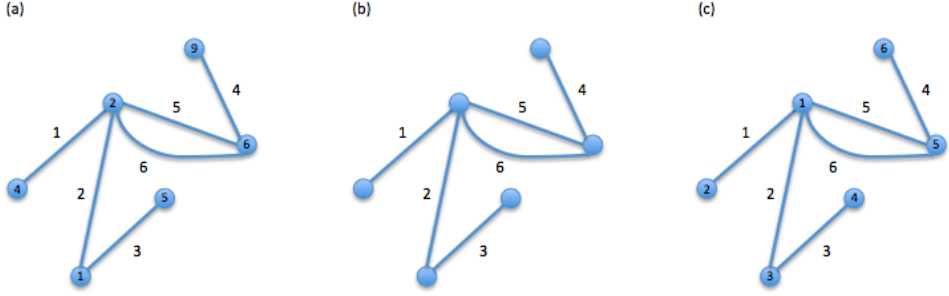

FIGURE 1. Illustration of the generic process of generating edge exchangeable networks in Section A. (a) Network representation of interaction process $X_1 = \{2, 4\}$, $X_2 = \{1, 2\}$, $X_3 = \{1, 5\}$, $X_4 = \{6, 9\}$, $X_5 = \{2, 6\}$, $X_6 = \{2, 6\}$. (b) Edge-labeled network $\mathcal{E}_X$ obtained by removing vertex labels from the vertex-edge-labeled network in Panel (a). (c) Vertex labeling of the network in Panel (b) according to its canonical selection function.

The $\operatorname{fin}_2(\mathbb{N})$-simplex consists of all $(f_{\{i,j\}})_{j \geq i \geq -1}$ such that $f_{\{i,j\}} \geq 0$ for all $j \geq i \geq -1$, $f_{\{-1,i\}} = 0$ for all $i \neq 0$, and $\sum_{j \geq i \geq -1} f_{\{i,j\}} = 1$. (The labels $-1$ and $0$ have a special status in the construction, as we shall see.) For any $f = (f_{\{i,j\}})_{j \geq i \geq -1}$ in the $\operatorname{fin}_2(\mathbb{N})$-simplex and $i \in \mathbb{N}$, we define

$$f_\bullet^{(i)} = \sum_{j=0}^\infty f_{\{i,j\}}$$

as the sum of masses involving element $i$.

Every $f = (f_{\{i,j\}})_{j \geq i \geq -1}$ in the $\operatorname{fin}_2(\mathbb{N})$-simplex determines a probability distribution on edge-labeled networks, denoted $\epsilon_f$, as follows. Let $X_1, X_2, \ldots$ be i.i.d. random pairs $\{i, j\}$ with

(1) $$\operatorname{pr}(X_k = \{i, j\} \mid f) = f_{\{i,j\}}, \quad j \geq i \geq -1.$$

Given $X_1, X_2, \ldots$, we define the selection function $\mathcal{X} : \mathbb{N} \to \operatorname{fin}_2(\mathbb{Z})$, where $\mathbb{Z} = \{\ldots, -1, 0, 1, \ldots\}$, as follows. We initialize with $m_0 = 0$. For $n \geq 1$, suppose $m_{n-1} = z \leq 0$. If $X_n$ contains no 0s, then we define $\mathcal{X}(n) = X_n$ and update $m_n = m_{n-1}$. If $X_n = \{0, j\}$ for some $j \geq 1$, then we put $\mathcal{X}(n) = \{z - 1, j\}$ and update $m_n = z - 1$. Otherwise, if $X_n = \{0, 0\}$, then we put $\mathcal{X}(n) = \{z - 1, z - 1\}$ and update $m_n = z - 1$; and if $X_n = \{-1, 0\}$, then we put $\mathcal{X}(n) = \{z - 1, z - 2\}$ and update $m_n = z - 2$. (Thus, events with $-1$ or $0$ involve vertices that appear once and never again. These are the "blips" that we previously ruled out in Crane and Dempsey (2016, Theorem 4.2).) We define $\mathcal{Y} = \mathcal{E}_\mathcal{X} \sim \epsilon_f$ to be the edge-labeled network induced by $\mathcal{X}$. See Figure 1 for an illustration of this procedure.

**Proposition A.1.** *The edge-labeled network $\mathcal{Y} = \mathcal{E}_X$ corresponding to $X_1, X_2, \ldots$ i.i.d. from* (1) *is edge exchangeable for all $f$ in the $\operatorname{fin}_2(\mathbb{N})$-simplex.*

For identifiability, we define the *ranked reordering of $f$* by $f^{\downarrow} = (f^{\downarrow}_{\{i,j\}})_{j \geq i \geq -1}$, the element of the $\operatorname{fin}_2(\mathbb{N})$-simplex obtained by putting $f^{\downarrow}_{\{-1,0\}} = f_{\{-1,0\}}$, $f^{\downarrow}_{\{0,0\}} = f_{\{0,0\}}$, and reordering elements $1, 2, \ldots$ so that $f^{(i)}_{\bullet} \geq f^{(i+1)}_{\bullet}$ for all $i \geq 1$ and then breaking ties $f^{(i)}_{\bullet} = f^{(i+1)}_{\bullet}$ by declaring that $(f_{\{i,i\}}, f_{\{i,i+2\}}, \ldots)$ comes before $(f_{\{i+1,i+1\}}, f_{\{i+1,i+2\}}, \ldots)$ in the lexicographic ordering. We write $\mathcal{F}^{\downarrow}$ to denote the space of rank reordered elements of the $\operatorname{fin}_2(\mathbb{N})$-simplex.

As the vertex labels other than $-1$ and $0$ are inconsequential, it is clear that $\epsilon_f$ and $\epsilon_{f'}$ determine the same distribution for any $f, f'$ for which $f^{\downarrow} = f'^{\downarrow}$. For any edge-labeled network $\mathcal{E}$, we write $|\mathcal{E}|^{\downarrow} \in \mathcal{F}^{\downarrow}$ to denote its *signature*, if it exists, as follows. Let $S_{\mathcal{E}} : \mathbb{N} \to \operatorname{fin}_2(\mathbb{N})$ be the canonical selection function for $\mathcal{E}$. For every $\{i, j\} \in \operatorname{fin}_2(\mathbb{N})$, $j \geq i \geq 1$, we define

$$f_{\{i,j\}}(\mathcal{E}) = \lim_{n \to \infty} n^{-1} \sum_{k=1}^{n} \mathbf{1}\{S_{\mathcal{E}}(k) = \{i, j\}\} \quad \text{and}$$

$$f^{(i)}_{\bullet}(\mathcal{E}) = \lim_{n \to \infty} n^{-1} \sum_{k=1}^{n} \mathbf{1}\{i \in S_{\mathcal{E}}(k)\},$$

if the limits exist. We also define

$$f_{\{0,i\}}(\mathcal{E}) = f^{(i)}_{\bullet} - \sum_{j=1}^{\infty} f_{\{i,j\}}(\mathcal{E}), \quad i \geq 1,$$

$$f_{\{0,0\}}(\mathcal{E}) = \lim_{n \to \infty} n^{-1} \sum_{k=1}^{n} \left( \sum_{\ell \geq 1} \mathbf{1}\{S_{\mathcal{E}}(k) = \{\ell, \ell\}\} \right) - \sum_{i=1}^{\infty} f_{\{i,i\}}(\mathcal{E}), \quad \text{and}$$

$$f_{\{-1,0\}}(\mathcal{E}) = \lim_{n \to \infty} n^{-1} \sum_{k=1}^{n} \left( \sum_{\ell, r \geq 1: \ell \neq r} \mathbf{1}\{S_{\mathcal{E}}(k) = \{\ell, r\}\} \right) - \sum_{j > i \geq 1} f_{\{i,j\}}(\mathcal{E}).$$

Provided each of the above limiting frequencies exists, we define $|\mathcal{E}| = (f_{\{i,j\}}(\mathcal{E}))_{j \geq i \geq -1}$ and $|\mathcal{E}|^{\downarrow} = (f_{\{i,j\}}(\mathcal{E}))^{\downarrow}_{j \geq i \geq -1}$.

(Note the role of $f_{\{0,0\}}(\mathcal{E})$ and $f_{\{-1,0\}}(\mathcal{E})$ for recording the residual proportion of loops and edges, respectively, that do not contribute to the limiting frequencies $f_{\{i,j\}}(\mathcal{E})$ for a given $j \geq i \geq 1$. For example, the interaction process $\mathcal{I} : \mathbb{N} \to \operatorname{fin}_2(\mathbb{N})$ given by $\mathcal{I}(i) = \{i, i\}$ for each $i \geq 1$ corresponds to the edge-labeled network $\mathcal{E} = \mathcal{E}_{\mathcal{I}}$ for which every edge is a loop at a distinct vertex. For any given $i \geq 1$, $f_{\{i,i\}}(\mathcal{E}) = 0$ and we have $f_{\{0,0\}}(\mathcal{E}) = 1$. Conversely, for $f = (f_{\{i,j\}})_{j \geq i \geq -1}$ with $f_{\{0,0\}} = 1$, the sequence $X_1, X_2, \ldots$ i.i.d. from (1) yields $\{0, 0\}, \{0, 0\}, \ldots$ with probability 1, whose associated interaction process has $\mathcal{X}(1) = \{-1, -1\}$, $\mathcal{X}(2) = \{-2, -2\}$, and $\mathcal{X}(n) = \{-n, -n\}$, for each $n \geq 1$, with associated edge-labeled network having all edges corresponding to a loop at a distinct vertex.)

**Theorem A.2.** *Let $\mathcal{Y} \in \mathfrak{C}_{\mathbb{N}}$ be an edge exchangeable network. Then there exists a unique probability measure $\phi$ on $\mathcal{F}^{\downarrow}$ such that $\mathcal{Y} \sim \epsilon_{\phi}$, where*

$$\epsilon_{\phi}(\cdot) = \int_{\mathcal{F}^{\downarrow}} \epsilon_f(\cdot) \phi(df). \tag{2}$$

*That is, every edge exchangeable network $\mathcal{Y}$ can be generated by first sampling $f \sim \phi$ and, given $f$, putting $\mathcal{Y} = \mathcal{E}_X$ for $X : \mathbb{N} \to \text{fin}_2(\mathbb{Z})$ constructed from $X_1, X_2, \ldots$ i.i.d. according to (1).*

Theorem 4.2 of Crane and Dempsey (2016) follows as a corollary to Theorem A.2 by ruling out blips, that is, by confining $\phi$ to the subset of $f \in \text{fin}_2(\mathbb{N})$ for which $f_{\{i,0\}} = 0$ for all $i \geq -1$.

A.1. **Proof of Theorem A.2.** We equip $\mathfrak{C}_{\mathbb{N}}$ with the product-discrete topology induced by the metric

$$d_{\mathfrak{C}_{\mathbb{N}}}(\mathcal{E}, \mathcal{E}') = 1/(1 + \sup\{n \in \mathbb{N} : \mathcal{E}|_{[n]} = \mathcal{E}'|_{[n]}\}), \quad \mathcal{E}, \mathcal{E}' \in \mathfrak{C}_{\mathbb{N}},$$

with convention $1/\infty = 0$, and $\mathcal{F}^{\downarrow}$ with the topology induced by

$$d_{\mathcal{F}^{\downarrow}}(f, f') = \sum_{j \geq i \geq -1} |f_{\{i,j\}} - f'_{\{i,j\}}|, \quad f, f' \in \mathcal{F}^{\downarrow}.$$

We then work with the respective Borel $\sigma$-fields induced by these topologies.

Let $\mathcal{Y} \in \mathfrak{C}_{\mathbb{N}}$ be an edge exchangeable random network, let $S_{\mathcal{Y}} : \mathbb{N} \to \text{fin}_2(\mathbb{N})$ be its canonical selection function, and let $\xi_1, \xi_2, \ldots$ be an i.i.d. sequence of Uniform$[0,1]$ random variables which are independent of $\mathcal{Y}$. Given $\mathcal{Y}$ and $(\xi_i)_{i \geq 1}$, we define $\mathcal{Z} : \mathbb{N} \to \text{fin}_2([0,1])$ by $\mathcal{Z}(n) = \{\xi_i, \xi_j\}$ on the event $S_{\mathcal{Y}}(n) = \{i, j\}$, for $n \geq 1$.

By independence of $\mathcal{Y}$ and $(\xi_i)_{i \geq 1}$ and edge exchangeability of $\mathcal{Y}$, $(\mathcal{Z}(n))_{n \geq 1}$ is an exchangeable sequence taking values in the Polish space $\text{fin}_2([0,1])$. By de Finetti's theorem, see, for example, Aldous (1985), there exists a unique measure $\mu$ on the space $\mathcal{P}(\text{fin}_2([0,1]))$ of probability measures on $\text{fin}_2([0,1])$ such that $\mathcal{Z} =_{\mathcal{D}} \mathcal{Z}^* = (\mathcal{Z}^*(n))_{n \geq 1}$ with

$$\text{pr}(\mathcal{Z}^* \in \cdot) = \int_{\mathcal{P}(\text{fin}_2([0,1]))} m^{\infty}(\cdot) \mu(dm),$$

where $m^{\infty}$ denotes the infinite product measure of $m$. In particular, there exists a random measure $\nu$ on $\text{fin}_2([0,1])$ such that

$$\text{pr}(\mathcal{Z} \in \cdot \mid \nu) = \nu^{\infty} \quad \text{a.s.}$$

Given $\nu$, we define

$$f_{\{i,j\}} = \nu(\{\{\xi_i, \xi_j\}\}), \quad i, j \geq 1,$$
$$f_\bullet^{(i)} = \nu(\{\{x, y\} \in \text{fin}_2([0,1]) : \xi_i \in \{x, y\}\}), \quad i \geq 1,$$
$$f_{\{0,i\}} = f_\bullet^{(i)} - \sum_{j=1}^\infty f_{\{i,j\}}, \quad i \geq 1,$$
$$f_{\{0,0\}} = \nu(\{\{u, u\} \in \text{fin}_2([0,1])\}) - \sum_{i=1}^\infty f_{\{i,i\}}, \quad \text{and}$$
$$f_{\{-1,0\}} = \nu(\{\{u, v\} \in \text{fin}_2([0,1]) : v \neq u\}) - \sum_{j>i\geq 1} f_{\{i,j\}}.$$

By construction $(f_{\{i,j\}})_{j\geq i\geq -1}$ is in the $\text{fin}_2(\mathbb{N})$-simplex and, therefore, $f^\downarrow \in \mathcal{F}^\downarrow$. (Note that $\mathcal{F}^\downarrow$ is a subset of the $\text{fin}_2(\mathbb{N})$-simplex and $f \mapsto f^\downarrow$ is measurable with respect to the Borel $\sigma$-field induced by the metric $d_{\mathcal{F}^\downarrow}(\cdot, \cdot)$ given above.)

Given $\nu$, we let $(\mathcal{Z}', S')$ be an i.i.d. copy of $(\mathcal{Z}, S_\mathcal{Y})$ and let $\mathcal{Y}' = \mathcal{E}_{S'}$ be the edge-labeled network induced by $S'$. We complete the proof by showing $\text{pr}(\mathcal{Y}' \in \cdot \mid \nu) = \epsilon_{f^\downarrow}$, for $f^\downarrow$ as defined above from $\nu$.

First, let $A = \{i \in \mathbb{N} : f_\bullet^{(i)} > 0\}$ and $\xi_A = \{\xi_i : i \in A\}$. It follows that

$$\text{pr}(\mathcal{Z}'(1) \cap \xi_A = \emptyset \mid \nu) = f_{\{0,0\}} + f_{\{-1,0\}} \quad \text{and} \quad \text{pr}(\mathcal{Z}'(1) \cap \xi_A = \{\xi_i\} \mid \nu) = f_{\{0,i\}}.$$

By exchangeability, $i \notin A$ implies $\xi_i$ appears at most once in $\mathcal{Z}$ with probability 1. We further have that

$$\text{pr}(\mathcal{Z}'_1 \cap \xi_A = \emptyset \text{ and } \mathcal{Z}'(1) = \{u, u\} \text{ for some } u \in [0,1] \mid \nu) = f_{\{0,0\}} \quad \text{and}$$
$$\text{pr}(\mathcal{Z}'(1) \cap \xi_A = \emptyset \text{ and } \mathcal{Z}'(1) = \{u, v\} \text{ for } u \neq v \mid \nu) = f_{\{-1,0\}}.$$

Now, define $\mathcal{X}' : \mathbb{N} \to \text{fin}_2(\mathbb{N} \cup \{-1, 0\})$ and the random selection function $S_{\mathcal{X}'} : \mathbb{N} \to \text{fin}_2(\mathbb{Z})$ as follows. Let $m_0 = 0$. For $n \geq 1$, suppose $m_{n-1} = z \leq 0$. If $\mathcal{Z}'(n) \cap \xi_A = \{\xi_i, \xi_j\}$ for some $i, j \in \mathbb{N}$, then put $\mathcal{X}(n) = S_X(n) = \{i, j\}$. If $\mathcal{Z}'(n) \cap \xi_A = \{\xi_i\}$ for some $i \geq 1$, then put $\mathcal{X}(n) = \{0, i\}$, $S_X(n) = \{z-1, i\}$, and $m_n = z - 1$. If $\mathcal{Z}'(n) \cap \xi_A = \emptyset$ and $\mathcal{Z}'(n) = \{u, u\}$ for some $u \in [0, 1]$, put $\mathcal{X}(n) = \{0, 0\}$, $S_\mathcal{X}(n) = \{z-1, z-1\}$, and $m_n = z - 1$. And if $\mathcal{Z}'(n) \cap \xi_A = \emptyset$ and $\mathcal{Z}'(n) = \{u, v\}$ for $u \neq v$, put $\mathcal{X}(n) = \{-1, 0\}$, $S_\mathcal{X}(n) = \{z-1, z-2\}$, and $m_n = z - 2$.

By construction, we have $S_{\mathcal{X}'} \equiv S'$ a.s. and, given $f$, $\mathcal{X}'$ is conditionally i.i.d. from distribution (1). The integral representation in (2) follows by de Finetti's theorem, completing the proof.

## Appendix B. Description of network datasets

Below we describe and provide references for some interaction datasets discussed in Crane and Dempsey (2016).

- Actors collaboration (Barabási and Albert, 1999): Network built from collaborations among actors in a given sample of movies. Each edge connects the actors who played a role in the corresponding movie.
- Enron email corpus (Klimt and Yang, 2004): Network built from a corpus of about 500,000 emails. Vertices are employees in the Enron Corporation with a directed edge from vertex $i$ to $j$ for each email sent from $i$ to $j$. Some versions of this dataset project edge multiplicities and ignore edge direction so that for each pair of vertices $i$ and $j$ there is an undirected edge between $i$ and $j$ if at least one email was exchanged between the two.
- Karate Club (Zachary, 1977): Network built from social interactions among 34 members of a karate club. Vertices are the members of the club and an edge between $i$ and $j$ corresponds to a social interaction between the two. The network exhibits no vertex sampling or growth since it is assumed all club members have been observed.
- Wikipedia voting (Leskovec et al., 2010): The Wikipedia voting network represents voting behavior for elections to the administrator role in Wikipedia. Vertices are Wikipedia users and a directed edge points from $i$ to $j$ if user $i$ voted for user $j$.
- US Airport (Colizza et al., 2007; Opsahl, 2011): The network from Opsahl (2011) is built from the flight map between all US airports in 2010. A directed edge from $i$ to $j$ indicates that a flight was scheduled from airport $i$ to airport $j$ in 2010. Edges are weighted by the number of seats on the scheduled flights. The network grows as a consequence of additional flights between airports.
- Co-authorship (Newman, 2001): Network built from co-authorship of preprints on the Condensed Matter section of arXiv between 1995 and 1999. Vertices are of two types, authors and papers, and edges only exist between vertices of a different type. An edge between $i$ (author) and $j$ (paper) indicates that $i$ is an author on paper $j$. The data is more succinctly represented by an edge-labeled network as in Crane and Dempsey (2016) by associating each article to an interaction involving all of its authors.
- UC Irvine (Opsahl and Panzarasa, 2009): Network built from UC Irvine online community. Vertices are active members of the community and a directed edge from $i$ to $j$ indicates that a message was sent from user $i$ to user $j$.
- Political blogs (Adamic and Glance, 2005): Network built from hyperlinks between political blogs. Vertices are websites (blogs) with a directed edge from $i$ to $j$ for every hyperlink from website $i$ to website $j$. Sampling is similar to the Facebook network.

Department of Statistics & Biostatistics, Rutgers University, 110 Frelinghuysen Avenue, Piscataway, NJ 08854, USA

*E-mail address*: hcrane@stat.rutgers.edu
*URL*: http://stat.rutgers.edu/home/hcrane

Department of Statistics, University of Michigan, 1085 S. University Ave, Ann Arbor, MI 48109, USA

*E-mail address*: wdem@umich.edu